\documentclass[twoside]{amsart}

\usepackage{amsmath,amssymb,amsfonts,amsthm,latexsym}
\usepackage[margin=1in]{geometry}
\usepackage{graphicx}
\usepackage[dvipsnames]{xcolor}
\usepackage{enumitem}
\usepackage{subcaption}
\usepackage{hyperref}
\usepackage{ytableau}
\usepackage{mathdots}
\usepackage{tikz-cd}
\usepackage{orcidlink}
\usepackage{mathtools}

\newtheorem{theorem}{Theorem}[section]
\newtheorem{lemma}[theorem]{Lemma}

\newtheorem{mainthm}{Theorem}

\theoremstyle{definition}

\theoremstyle{remark}

\numberwithin{equation}{section}

\newcommand{\C}{\mathbb{C}}
\newcommand{\Sp}{\mathbb{S}}

\newcommand{\R}{\mathbb{R}}

\newcommand{\Z}{\mathbb{Z}}

\newcommand{\Real}{\operatorname{Re}}
\newcommand{\Imaginary}{\operatorname{Im}}

\newcommand{\U}{\operatorname{U}}

\newcommand{\diag}{\operatorname{diag}}

\newcommand{\OO}{\operatorname{O}}
\newcommand{\SO}{\operatorname{SO}}

\newcommand{\vol}{\text{vol}}
\newcommand{\Vol}{\operatorname{Vol}}
\newcommand{\sign}{\operatorname{sign}}

\DeclarePairedDelimiter\abs{\lvert}{\rvert}%
\DeclarePairedDelimiter\norm{\lVert}{\rVert}%
\makeatletter
\let\oldabs\abs
\def\abs{\@ifstar{\oldabs}{\oldabs*}}
\let\oldnorm\norm
\def\norm{\@ifstar{\oldnorm}{\oldnorm*}}
\makeatother

\begin{document}

\title[Hilbert measures of non-coregular representations]{Hilbert measures of non-coregular representations of the circle and orthogonal groups}

\author[L.~Liebow]{Luca Liebow}
\address{Department of Mathematics and Statistics,
Skidmore College,
815 North Broadway,
Saratoga Springs, NY 12866, USA}

\author[B.~Roberge]{Belle Roberge}
\address{Department of Mathematics and Statistics,
Skidmore College,
815 North Broadway,
Saratoga Springs, NY 12866, USA}

\author[C.~Seaton]{Christopher Seaton\,\orcidlink{0000-0001-8545-2599}}
\address{Department of Mathematics and Statistics,
Skidmore College,
815 North Broadway,
Saratoga Springs, NY 12866, USA}
\email{cseaton@skidmore.edu}


\thanks{All authors were supported by the Skidmore College Faculty Student Summer Research Program.
C.S. was supported by an AMS-Simons Research Enhancement Grant for PUI Faculty.}

\begin{abstract}
We give explicit constructions of the Hilbert measure, which transports the integral of an invariant smooth function from a representation to the image of the Hilbert embedding, for certain representations of the circle, orthogonal, and special orthogonal groups that are non-coregular, meaning that the algebra of invariants is not a polynomial algebra. These include all $4$-dimension real circle representations, a family of $6$-dimensional circle representations, and $k$ copies of the defining representations of the orthogonal and special orthogonal groups when $k$ is large enough that the representation is not coregular. For the circle representations, we also construct the Hilbert measure on the singular symplectic quotient at level zero.
\end{abstract}


\setcounter{tocdepth}{3}
\maketitle

\tableofcontents


\section{Introduction}
\label{sec:intro}

Let $G$ be a compact Lie group and let $V$ be a finite-dimensional real $G$-representation.
Given a minimal generating set $u_1,\ldots,u_\ell$ for the real $G$-invariant polynomials on $V$, which we refer to as a \emph{Hilbert basis},
the corresponding \emph{Hilbert embedding} is the map $H\colon V\to\R^\ell$ given by $H = (u_1,\ldots,u_\ell)$. The Hilbert embedding induces
a homeomorphism of the orbit space $G\backslash V$ onto its image $X = H(V)$ and hence yields a model for the orbit space as a semialgebraic
subset of Euclidean space. See \cite{BierstoneBook,ProcesiSchwarz} for more details.

The semialgebraic set $X$ is not merely a model for the orbit space as a topological space. By \cite[Thrm.~A]{BierstoneLiftingIsotopies},
see also \cite[Thrm.~2.5]{BierstoneBook}, the stratification of the orbit space $G\backslash V$ by orbit types coincides with the
minimum primary stratification of $X$ as a semialgebraic, hence semianalytic subset of $\R^\ell$ constructed in
\cite[Sec.~25]{LojasiewiczES-A}; see also \cite[Thrm.~4.9]{MatherStratMappings}.
A crucial result for this paper is \cite[Thrm.~1]{GWSdiff}, that the algebra of smooth functions on
$X$ given by the restrictions of smooth functions on $\R^\ell$ coincides, via $H$, with the algebra of smooth functions on $G\backslash V$ given
by functions that pull back to smooth functions on $V$; see also \cite[Thrm.~1]{MatherDiffInvariants}. In particular, any smooth $G$-invariant
function on $V$ can be written as a smooth function of the Hilbert basis $u_1,\ldots,u_\ell$.

Now, let $f$ be a smooth $G$-invariant real-valued function on $V$ and suppose that $f$ is integrable on $V$. The result of Schwarz described
above implies that we can express $f = F(u_1,\ldots,u_\ell)$ for a smooth function $F$ on $\R^\ell$ and hence consider $f$ as a function on $X$.
This suggests that we can express the integral
\begin{equation}
\label{eq:MainIntegral}
    \int_V f \, d\vol_V
    =   \int_X F \, \Lambda
\end{equation}
where $\Lambda$ is a differential form on an open dense subset of $X$. We refer to $\Lambda$ as the \emph{Hilbert measure} associated to
the Hilbert embedding $H$. Note that $\Lambda$ is oriented by a choice of orientation on $V$.

In \cite{HerbigSeatonWhitesell}, the Hilbert measure was computed for $k$ copies of the defining representation of the real orthogonal group
$\OO_m$ when $k\leq m$ (as well as the corresponding restrictions to $\SO_m$ when $k\leq m-1$). These representations are \emph{coregular},
meaning that the algebra of invariants is a polynomial algebra or, equivalently, the elements of the Hilbert basis are algebraically independent.
It follows that the image $X$ of the Hilbert embedding is a region in $\R^\ell$ of dimension $\ell$, i.e., restricted only by the inequalities determined
by \cite[Thrm.~0.10]{ProcesiSchwarz} and not cut down by inequalities. In the coregular case, the Hilbert measure is of the form
\[
    \Lambda
    =   \lambda \, du_1 \wedge \cdots \wedge du_\ell
\]
for a function $\lambda$ on $X$, and hence is determined by the function $\lambda$.

Roughly speaking, ``most" representations of compact Lie groups are not coregular; see
\cite{AdamovichGolovina,KacPopovVinbergLibre,LittelmannKoregulaere,GWSSimpleRegularRingsInvars}.
In the non-coregular case, $X$ is cut down by equations and inequalities, so the Hilbert measure cannot be described by a function $\lambda$.
The purpose of this paper is to extend the techniques of \cite{HerbigSeatonWhitesell} to non-coregular examples and explicitly compute
the differential form $\Lambda$. The resulting differential form is well-defined on the open dense subset of $X$ given by the principal
orbit type, the top stratum in the stratification by orbit types, and hence a smooth manifold. The complement of this subset has measure zero,
so is sufficient to define the integral \eqref{eq:MainIntegral}. Singularities of the Hilbert measure along lower-dimensional strata were
already observed in the coregular case; see \cite[pp.~3]{HerbigSeatonWhitesell}.

The outline and main results of this paper are as follows. In Section~\ref{sec:S1-ab},
we consider the $4$-dimensional real representations of the circle $\Sp^1$. As explained at the beginning of that section,
assuming such a representation has no trivial factors, it is isomorphic to the underlying real representation of a unitary representation
with weight vector $(-a,b)$. After describing the image $X$ of the Hilbert embedding in Section~\ref{subsec:S1-ab_ImageHilb}, we prove the following
in Sections~\ref{subsec:S1-ab_COV1} and \ref{subsec:S1-ab_COV2}.

\begin{mainthm}
\label{mainthm:S1-ab}
Let $V\simeq\C^2$ be the $\Sp^1$-representation with weight vector $(-a,b)\in\Z^2$ where $a$ and $b$ are positive integers and $\gcd(a,b) = 1$.
Let $H\colon V\to\R^4$ be the Hilbert embedding defined using the Hilbert basis \eqref{eq:S1-ab_RInvars} and let $X = H(V)$ be the image
of the Hilbert embedding. The Hilbert measure on $X$ associated to $V$ and $H$ is
\begin{equation}
\label{eq:Main_S1-ab}
\begin{aligned}
    \Lambda_X
    &=
    - \frac{\pi }
        {2 u_4} \, du_1\wedge du_2 \wedge du_3
    &&=
    \frac{\pi}
        {2 u_3} \, du_1\wedge du_2 \wedge du_4
    \\&=
    \frac{\pi }
        {a u_1^b u_2^{a-1}} \, du_1\wedge du_3 \wedge du_4
    &&=
    - \frac{\pi}
        {b u_1^{b - 1} u_2^a} \, du_2\wedge du_3 \wedge du_4.
\end{aligned}
\end{equation}
\end{mainthm}

In Section~\ref{subsec:S1-ab_Symp}, we consider the corresponding symplectic quotient and prove the following.

\begin{mainthm}
\label{mainthm:S1-ab_Symp}
Let $V\simeq\C^2$ be the $\Sp^1$-representation with weight vector $(-a,b)\in\Z^2$ where $a$ and $b$ are positive integers.
Let $Z$ be the zero fiber of the homogeneous quadratic moment map $J\colon V\to\mathfrak{t}^\ast$,
let $H_Z\colon Z\to\R^3$ be the Hilbert embedding described in Section~\ref{subsubsec:S1-ab_SympImageHilbDS}, and let $X_Z = H(Z)$ be the image
of the Hilbert embedding. The Hilbert measure on the singular symplectic quotient associated to $H_Z$ is
\begin{equation}
\label{eq:Main_S1-ab_Symp}
    \Lambda_{X_Z}
    =
    - \frac{ a \pi \sqrt{a(a + b) u_1 }  }
        { b u_4 }
        \, du_1\wedge du_3
    =
    \frac{ a \pi \sqrt{a(a + b)u_1}  }
        { b u_3}
        \, du_1\wedge du_4
    =
    \frac{ 2\pi b^{a-1}}
        { a^{a-3/2} \sqrt{a + b} u_1^{a + b - 3/2} }
        \, du_3\wedge du_4.
\end{equation}
\end{mainthm}

In Section~\ref{sec:S1-1aa}, we consider a family of $6$-dimensional real representations of $\Sp^1$, those given by weight vectors $(-1,a,a)$
for a positive integer $a$. We once again describe the image of the Hilbert embedding in Section~\ref{subsec:S1-1aa_ImageHilb},
and in Sections~\ref{subsec:S1-1aa_COV1} and \ref{subsec:S1-1aa_COV2}, we prove the following.

\begin{mainthm}
\label{mainthm:S1-1aa}
Let $V\simeq\C^3$ be the $\Sp^1$-representation with weight vector $(-1,a,a)\in\Z^3$ where $a$ is a positive integer.
Let $H\colon V\to\R^9$ be the Hilbert embedding defined using the Hilbert basis \eqref{eq:S1-1aa_RInvars} and let $X = H(V)$ be the image
of the Hilbert embedding. The Hilbert measure on $X$ associated to $V$ and $H$ is
\begin{equation}
\label{eq:Main_S1-1aa}
     \Lambda_X
     =
     \frac{\pi}{ 4 u_5 u_6 } \, du_1\wedge du_2 \wedge du_3 \wedge du_4 \wedge du_7
     =
     \frac{\pi}{ 4 u_4 u_7 } \, du_1\wedge du_2 \wedge du_3 \wedge du_5 \wedge du_6.
\end{equation}
\end{mainthm}

We also consider the symplectic quotient in Section~\ref{subsec:S1-1aa_Symp}, proving the following.

\begin{mainthm}
\label{mainthm:S1-1aa_Symp}
Let $V\simeq\C^3$ be the $\Sp^1$-representation with weight vector $(-1,a,a)\in\Z^3$ where $a$ is a positive integer.
Let $Z$ be the zero fiber of the homogeneous quadratic moment map $J\colon V\to\mathfrak{t}^\ast$,
let $H_Z\colon Z\to\R^8$ be the Hilbert embedding described in Section~\ref{subsec:S1-1aa_Symp}, and let $X_Z = H(Z)$ be the image
of the Hilbert embedding. The Hilbert measure on the singular symplectic quotient associated to $H_Z$ is
\begin{equation}
\label{eq:Main_S1-1aa_Symp}
    \Lambda_{X_Z}
    =
    - \frac{\pi \sqrt{ (1 + a) u_1 } }
        { 2 a^2 u_4 u_7 }
        \, du_1\wedge du_2 \wedge du_5 \wedge du_6        
    =
    - \frac{\pi \sqrt{ (1 + a) u_1 } }
        { a u_5 (u_1 - 2a u_2) }
        \, du_1\wedge du_4 \wedge du_8 \wedge du_9.
\end{equation}
\end{mainthm}

In Section~\ref{sec:Om_kRm} we turn to the case of $k$ copies of the defining representation of the real orthogonal group $\OO_m$ with $k > m$.
Recall that the volume of the $(n-1)$-dimensional sphere $S^{n-1}$ is given by $\Vol(S^{n-1}) = \frac{2\pi^{n/2}}{\Gamma(n/2)}$.
Following the outline of the previous sections, we describe the image of the Hilbert embedding in Section~\ref{subsec:Om_ImageHilb}
and compute the Hilbert measure in Section~\ref{subsec:Om_COV1} and \ref{subsec:Om_COV2}. In Section~\ref{sec:SOm_kRm}, we describe the modifications
for the restriction of the representation to $\SO_m$, yielding the following.

\begin{mainthm}
\label{mainthm:Om-SOm}
Let $V_k \simeq\R^{km}$ be the $\OO_m$-representation given by $k$ copies of the defining representation with $k > m$.
Let $H\colon V_k\to\R^{{k+1\choose 2}}$ be the Hilbert embedding defined in \eqref{eq:Om-HilbertEmbed}
and let $X = H(V)$ be the image of the Hilbert embedding. Ordering the variables $u_{i,j}$ lexicographically,
the Hilbert measure associated to $V$ and $H$ is the form
\begin{equation}
\label{eq:Main_Om}
\begin{aligned}
    \Lambda_X
    &=
    \frac{ \pi^{(m-1)(m+2)/4}
            \prod\limits_{i=1}^{m-2} \left( \frac{ \Gamma\left(\frac{i+1}{2}\right) }{ \Gamma\left(\frac{i+2}{2}\right)} \right)^{m-i-1} }
        {(\det G_m)^{(r+1)/2}}
        \, \bigwedge\limits_{\substack{1\leq i\leq k \\ 1\leq j \leq {\min\{i,m\}} }}  du_{i,j}
    \\&=
    \frac{ \prod\limits_{i=1}^m \Vol(S^{i-1}) }{2^m (\det G_m)^{(r+1)/2}}
        \, \bigwedge\limits_{\substack{1\leq i\leq k \\ 1\leq j \leq {\min\{i,m\}} }}  du_{i,j}.
\end{aligned}
\end{equation}

Restricting to the $\SO_m$-representation, using the Hilbert embedding $H$ in \eqref{eq:SOm-HilbertEmbed}, and letting $X$ denote the image of $H$,
the Hilbert measure is given by
\begin{equation}
\label{eq:Main_SOm}
\begin{split}
    \Lambda_X &=
    \frac{ \prod\limits_{i=2}^m \Vol(S^{i-1}) }{\sign(\det\nolimits_{[m]})^{r+1} 2^m (\det G_m)^{(r+1)/2}}
        \, \bigwedge\limits_{\substack{1\leq i\leq k \\ 1\leq j \leq {\min\{i,m\}} }}  du_{i,j}
    \\&=
    \frac{ \pi^{(m-1)(m+2)/4}
        \prod\limits_{i=1}^{m-2} \left( \frac{ \Gamma\left(\frac{i+1}{2}\right) }{ \Gamma\left(\frac{i+2}{2}\right)} \right)^{m-i-1} }
        {2\sign(\det\nolimits_{[m]})^{r+1} (\det G_m)^{(r+1)/2}}
        \, \bigwedge\limits_{\substack{1\leq i\leq k \\ 1\leq j \leq {\min\{i,m\}} }}  du_{i,j}.
\end{split}
\end{equation}
\end{mainthm}

The coefficients in Equations~\eqref{eq:Main_Om} and \eqref{eq:Main_SOm} can also be written in terms of the volume of $\OO_m$ as,
by \cite[Thrm.~2.24]{ZhangVolumes},
\[
    \Vol(\OO_m)
        =   \prod\limits_{i=1}^m \Vol(S^{i-1})
        =   \frac{2^m \sqrt{\pi}^{m+1\choose 2} }{\prod_{j=1}^m \Gamma(j/2) }.
\]

The outline of each computation is similar. We identify a fundamental domain and use it to parametrize an open dense subset of $V$ or $Z$.
Because we are interested in integrals and hence can ignore sets of measure zero, we use the term \emph{fundamental domain} rather loosely;
a more correct term would be \emph{fundamental domain for an open dense subset}.
This is referred to as the ``first change of variables" and allows us to integrate out the parameters in the direction of the orbit.
The result is a simplified integral on this fundamental domain; see
Equations~\eqref{eq:S1-ab_Int1}, \eqref{eq:S1-ab_SympInt1}, \eqref{eq:S1-1aa_Int1}, \eqref{eq:S1-1aa_SympInt1},
\eqref{eq:Om_Int1}, and \eqref{eq:SOm_Int1}. The ``second change of variables" then transforms the result to an integral on the
image of the Hilbert embedding, yielding the Hilbert measure $\Lambda$.


\section*{Acknowledgements}

This paper developed from L.L and B.R.'s research projects conducted through the Skidmore College Faculty Student Summer Research Program
and Skidmore Mathematics and Statistics Department, and the authors gratefully acknowledge the support of the department and college for these activities.
All authors were supported by the Skidmore College Faculty Student Summer Research Program, and
C.S. was supported by an AMS-Simons Research Enhancement Grant for PUI Faculty.
C.S. would like to honor the memory of Hans-Christian Herbig, who proposed and with whom he initiated the line of inquiry continued here.


\section{Circle representations on $\C^2$}
\label{sec:S1-ab}

When $G = \Sp^1$ is the group of complex units, a unitary representation $V$ of $G$ is described in a \emph{weight vector} $(a_1,\ldots,a_n)\in\Z^n$.
Given such a weight vector, there are corresponding coordinates $(z_1,\ldots,z_n)$ for $V\simeq\C^n$ such that
the action is given by the homomorphism mapping $t\in\Sp^1$ to the diagonal matrix $\diag(t^{a_1}, \ldots, t^{a_n}) \in\U(n,\C)$.
See \cite[Ch.~II, Sec.~8]{BroeckerTomDieck} for more details.

In this section, we consider the underlying real representation of the unitary representation with weight vector $(-a, b)$ where $a,b\in\Z$ are positive.
We assume that $\gcd(a,b) = 1$; otherwise, the finite subgroup of $\Sp^1$ of order $\gcd(a,b)$ acts trivially on $V$
and we may reduce to the quotient of $\Sp^1$ by this subgroup, an isomorphic group that acts with weight vector $(-a/\gcd(a,b), b/\gcd(a,b))$,
with no changes to our results. We prove Theorem~\ref{mainthm:S1-ab} in Sections~\ref{subsec:S1-ab_COV1} and \ref{subsec:S1-ab_COV2}.

Though the unitary representations associated to positive and negative weights are distinct, the underlying real representations with weight vectors
$(a)$ and $(-a)$ are isomorphic; see \cite[Ch.~II, Prop.~9.5]{BroeckerTomDieck}. Hence, any $4$-dimensional real representation of the circle
with no trivial factors is isomorphic to a representation with weight vector $(-a, b)$.

In Section~\ref{subsec:S1-ab_Symp}, we will consider the singular symplectic quotient at level zero associated to this representation and prove
Theorem~\ref{mainthm:S1-ab_Symp}.
For a unitary representation $V$ with weight vector $(a_1,\ldots,a_n)$, letting $\mathfrak{t}$ denote the Lie algebra of $\Sp^1$,
the \emph{moment map} $J\colon V\to\mathfrak{t}^\ast\simeq\R$ is given by
\[
    J(z_1,\ldots,z_n) = - \frac{1}{2}(a_1\abs{z}_1 + \cdots a_n\abs{z}_n).
\]
The \emph{shell} is the subset $J^{-1}(0)$, which is $\Sp^1$-invariant, and
the \emph{symplectic quotient at level zero} is given by $\Sp^1\backslash J^{-1}(0)$.
See \cite[Sec.~4.4]{GuilleminGinzburgKarshon} and \cite{SjamaarLerman} for more details.
The symplectic quotient depends on the isomorphism class of the unitary representation, not the underlying real representation, and hence
depends on the sign of the weights. However, if both weights have the same sign, then $J^{-1}(0)$ is a point, so that the nontrivial cases
associated to $4$-dimensional real representation of the circle are isomorphic to those considered here.


\subsection{The image of the Hilbert embedding}
\label{subsec:S1-ab_ImageHilb}

In complex coordinates $(z_1, z_2)$ for $V$, a minimal set of generating invariants consists of
\begin{equation}
\label{eq:S1-ab_CInvars}
\begin{split}
    m_1 &=      z_1 \overline{z_1},
    \\
    m_2 &=      z_2 \overline{z_2},
    \\
    p &=      z_1^b z_2^a,
    \\
    q &=      \overline{z_1}^b \overline{z_2}^a,
\end{split}
\end{equation}
and they satisfy the single relation
\begin{equation}
\label{eq:S1-ab_CRelat}
    m_1^b m_2^a -   p q = 0.
\end{equation}
In real coordinates $(\mathbf{x}, \mathbf{y}) = (x_1, y_1, x_2, y_2)$ for $V$ where $z_j = x_j + i y_j$ for $j = 1,2$,
the invariants $m_j$ are real and given by $m_j = x_j^2 + y_j^2$.
The remaining invariants are given by the real and imaginary parts of $p$, (which coincide with the real and imaginary parts of
$q$ up to a sign).
Hence the real Hilbert basis is given by
\begin{equation}
\label{eq:S1-ab_RInvars}
\begin{aligned}
    u_1 &=      x_1^2 + y_1^2,
    \qquad\qquad
    u_2 =       x_2^2 + y_2^2,
    \\
    u_3 &=      \Real p
        =       \frac{1}{2}\big( (x_1 + i y_1)^b + (x_1 - i y_1)^b \big)\big( (x_2 + i y_2)^a + (x_2 - i y_2)^a \big),
    \\
    u_4 &=      \Imaginary p
        =       \frac{1}{2i}\big( (x_1 + i y_1)^b - (x_1 - i y_1)^b \big)\big( (x_2 + i y_2)^a - (x_2 - i y_2)^a \big).
\end{aligned}
\end{equation}
It follows that the Hilbert embedding $H\colon V\to\R^4$ is given by
$H(\mathbf{x}, \mathbf{y}) = (u_1, u_2, u_3, u_4)$. From Equation~\eqref{eq:S1-ab_CRelat}, one computes that the real invariants satisfy the single relation
\begin{equation}
\label{eq:S1-ab_RRelat}
    u_1^b u_2^a - u_3^2 - u_4^2 = 0.
\end{equation}

\begin{lemma}
\label{lem:S1-ab_ImHilb}
The image of the Hilbert embedding associated to the invariants \eqref{eq:S1-ab_RInvars} is
\begin{equation}
\label{eq:S1-ab_ImHilb}
    X = \{
    (u_1, u_2, u_3, u_4)\in\R^4 :
    u_1^b u_2^a - u_3^2 - u_4^2 = 0,
    \;
    u_1 \geq 0,
    \;
    u_2 \geq 0
    \}.
\end{equation}
\end{lemma}
\begin{proof}
By \cite[Thrm~0.10 and Prop.~0.6]{ProcesiSchwarz}, the inequalities that define $X$ as a subset of the variety corresponding to the relation
\eqref{eq:S1-ab_RRelat} correspond to the non-negativity of the principal minors of the matrix
\[
    (\nabla u_j \cdot \nabla u_k)
    =
    \begin{pmatrix}
        4 u_1 & 0 & 2 b u_3 & 2 b u_4 \\
        0 & 4 u_2 & 2 a u_3 & 2 a u_4 \\
        2 b u_3 & 2 a u_3 & u_2^{a-1} u_1^{b-1}
        (a^2 u_1+b^2 u_2) & 0 \\
        2 b u_4 & 2 a u_4 & 0 & u_2^{a-1} u_1^{b-1} (a^2 u_1+b^2 u_2)
    \end{pmatrix}.
\]
Computing the principal minors of this matrix and simplifying the resulting inequalities, they can be reduced to
\[
    u_1 \geq 0,
    \quad
    u_2 \geq 0,
    \quad
    u_1^b u_2^a - u_3^2 \geq 0,
    \quad
    u_1^b u_2^a - u_3^2 - u_4^2 \geq 0.
\]
The last inequality is always satisfied due to the relation \eqref{eq:S1-ab_RRelat}, which also implies
$u_1^b u_2^a - u_3^2 = u_4^2 \geq 0$, so on the set of points satisfying \eqref{eq:S1-ab_RRelat}, the only non-redundant
inequalities are $u_1 \geq 0$ and $u_2 \geq 0$.
\end{proof}


\subsection{The first change of variables}
\label{subsec:S1-ab_COV1}

Let $f(\mathbf{x}, \mathbf{y})$ be a smooth, real-valued, integrable function on $V$ and assume $f$ is invariant under the action of $\Sp^1$. Our goal is to
express the integral
\begin{equation}
\label{eq:S1-ab_Int0}
    \int_V f(\mathbf{x}, \mathbf{y}) \, dx_1 \wedge dy_1 \wedge dx_2\wedge dy_2
\end{equation}
as an integral over the image of the Hilbert embedding. First, we restrict to a fundamental domain (for an open dense subset) of $V$ as follows.

Let $(\mathbf{x}, \mathbf{y}) \in V$ with $x_j + i y_j \neq 0$ for $j = 1,2$, and then, recalling that $\gcd(a, b) = 1$, there is a unique
element of the form $(w_1, 0, w_2, w_3)$ of the $\Sp^1$-orbit of $(\mathbf{x}, \mathbf{y})$ such that $w_1 = \sqrt{x_1^2 + y_1^2} > 0$ and
the argument of $w_2 + i w_3$ is contained in $[0, 2\pi/a)$. Choosing the unique $\theta\in(-\pi,\pi]$ such that
$e^{i\theta}\ast(w_1, 0, w_2, w_3) = (\mathbf{x}, \mathbf{y})$, we have
\begin{equation}
\label{eq:S1-ab_COV1}
    (\mathbf{x}, \mathbf{y})
    =
    \big(   w_1 \cos(a\theta),
            - w_1 \sin(a\theta),
            w_2 \cos(b\theta) - w_3 \sin(b\theta),
            w_2 \sin(b\theta) + w_3 \cos(b\theta)
    \big).
\end{equation}
Noting that $\theta$ moves in the negative direction so that $(\theta, w_1, w_2, w_3)$ has the same orientation as $(\mathbf{x}, \mathbf{y})$,
the Jacobian of the change of variables is
\begin{equation}
\label{eq:S1-ab_COV1Jacobian}
    \frac{\partial (\mathbf{x}, \mathbf{y})}{ \partial (\theta, w_1, w_2, w_3)}
    =
    \begin{pmatrix}
        - a w_1\sin(a\theta)    &   \cos(a\theta)     &   0   &   0
        \\
        - a w_1\cos(a\theta)    &   - \sin(a\theta)     &   0   &   0
        \\
        - b w_2 \sin(b\theta) - b w_3 \cos(b\theta)   & 0   &   \cos(b\theta)   &   -\sin(b\theta)
        \\
        b w_2 \cos(b\theta) - b w_3 \sin(b\theta)     & 0   &   \sin(b\theta)   &   \cos(b\theta)
    \end{pmatrix}
\end{equation}
and has determinant $a w_1$ so that
\[
    dx_1\wedge dy_1 \wedge dx_2 \wedge dy_2
    =
    a w_1 d\theta \wedge dw_1 \wedge dw_2 \wedge dw_3.
\]
Composing $f$ with the change of variables \eqref{eq:S1-ab_COV1} yields a smooth function $g$ of $(\theta, w_1, w_2, w_3)$,
and $\Sp^1$-invariance implies that $g$ is independent of $\theta$. Hence, we have $g(w_1, w_2, w_3) = f(\mathbf{x}, \mathbf{y})$, and letting
\[
    W = \{ (w_1, 0, w_2, w_3) \in V : w_1 > 0, \arg (w_2 + i w_3) \in [0,2\pi/a) \},
\]
we can express the integral in Equation~\eqref{eq:S1-ab_Int0} as
\begin{align}
    \nonumber
    &\int_{(-\pi,\pi]\times W} g(w_1, w_2, w_3) a w_1 \, d\theta \wedge dw_1 \wedge dw_2 \wedge dw_3
    \\
    \nonumber
    &\qquad=    \int_{-\pi}^{\pi} d\theta \int_W g(w_1, w_2, w_3) a w_1 \, dw_1 \wedge dw_2 \wedge dw_3
    \\
    \label{eq:S1-ab_Int1}
    &\qquad=    2a\pi \int_W g(w_1, w_2, w_3) w_1 \, dw_1 \wedge dw_2 \wedge dw_3.
\end{align}
This gives an expression for the original integral \eqref{eq:S1-ab_Int0} on the fundamental domain $W$.


\subsection{The second change of variables}
\label{subsec:S1-ab_COV2}

The next step is to express the integral \eqref{eq:S1-ab_Int1} as an integral on the image of the Hilbert embedding $X$ given in
Equation~\eqref{eq:S1-ab_ImHilb}. Considering the $w_j$ as parameters for $X$, the parametrization is given by substituting
the $x_1 = w_1$, $y_1 = 0$, $x_2 = w_2$, and $y_2 = w_3$ into the definition of the $u_k$ given in Equation~\eqref{eq:S1-ab_RInvars},
yielding
\begin{equation}
\label{eq:S1-ab_RInvarsInW}
\begin{aligned}
    u_1 &=      w_1^2
    &
    u_2 &=      w_2^2 + w_3^2,
    \\
    u_3 &=      \frac{w_1^b}{2}\big( (w_2 + i w_3)^a + (w_2 - i w_3)^a \big),
    \qquad&
    u_4 &=      \frac{w_1^b}{2i}\big( (w_2 + i w_3)^a - (w_2 - i w_3)^a \big).
\end{aligned}
\end{equation}
To compute the volume form on $X$ in terms of this parametrization, the Jacobian
$\frac{\partial\mathbf{u}}{\partial\mathbf{w}} = \frac{\partial (u_1, u_2, u_3, u_4)}{ \partial (w_1, w_2, w_3)}$ is given by
\[
    \frac{\partial\mathbf{u}}{\partial\mathbf{w}}
    =
    \begin{pmatrix}
        2w_1 & 0 & \frac{bw_1^{b-1}}{2}\big( (w_2 + i w_3)^a + (w_2 - i w_3)^a \big) & \frac{b w_1^{b-1}}{2i}\big( (w_2 + i w_3)^a - (w_2 - i w_3)^a \big)
        \\
        0 & 2w_2 & \frac{a w_1^b}{2}\big( (w_2 + i w_3)^{a-1} + (w_2 - i w_3)^{a-1} \big) & \frac{a w_1^b}{2i}\big( (w_2 + i w_3)^{a-1} - (w_2 - i w_3)^{a-1} \big)
        \\
        0 & 2w_3 & \frac{a i w_1^b}{2}\big( (w_2 + i w_3)^{a-1} - (w_2 - i w_3)^{a-1} \big) & \frac{a w_1^b}{2}\big( (w_2 + i w_3)^{a-1} + (w_2 - i w_3)^{a-1} \big)
    \end{pmatrix}^T,
\]
from which one computes that
\begin{align*}
    \det \Bigg(\Big(\frac{\partial\mathbf{u}}{\partial\mathbf{w}}\Big)^T \frac{\partial\mathbf{u}}{\partial\mathbf{w}}\Bigg)
    &=
    4 a^2 w_1^{2b} (w_2^2 + w_3^2)^{a-1}\Big(
        4 w_1^2 (w_2^2 + w_3^2) + a^2 w_1^{2b+2} (w_2^2 + w_3^2)^{a-1} + b^2 w_1^{2b-2} (w_2^2 + w_3^2)^{a+1} \big) \Big),
    \\&=
    4 a^2 u_1^b u_2^{a-1}\big(
        4 u_1 u_2 + a^2 u_1^{b+1} u_2^{a-1} + b^2 u_1^{b-1} u_2^{a+1} \big).
\end{align*}
It follows that the volume element on $X$ in terms of this parametrization is given by
\[
    d\vol_X  =
    2 a u_1^{b/2} u_2^{(a-1)/2}\sqrt{
        4 u_1 u_2 + a^2 u_1^{b+1} u_2^{a-1} + b^2 u_1^{b-1} u_2^{a+1}
        } \, dw_1 \wedge dw_2 \wedge dw_3.
\]
Let $F$ be a smooth function such that $F(u_1, u_2, u_3, u_4) = g(w_1, w_2, w_3) = f(\mathbf{x}, \mathbf{y})$;
the existence of such a function is guaranteed by \cite[Thrm.~1]{GWSdiff}. Then
\[
    2a\pi \int_W g(w_1, w_2, w_3) w_1 \, dw_1 \wedge dw_2 \wedge dw_3
    =
    \int_X F(u_1, u_2, u_3, u_4) \kappa \, d\vol_X
\]
implies that
\begin{align*}
    \kappa
    &=
    \frac{2 a \pi w_1}
        {2 a u_1^{b/2} u_2^{(a-1)/2} \sqrt{
            4 u_1 u_2 + a^2 u_1^{b+1} u_2^{a-1} + b^2 u_1^{b-1} u_2^{a+1}
        } }
    \\&=
    \frac{\pi}
        {u_1^{(b-1)/2} u_2^{(a-1)/2} \sqrt{
            4 u_1 u_2 + a^2 u_1^{b+1} u_2^{a-1} + b^2 u_1^{b-1} u_2^{a+1}
        } }.
\end{align*}
Therefore, we can express the integral \eqref{eq:S1-ab_Int1} as an integral on $X$ as
\begin{equation}
\label{eq:S1-ab_Int2}
    \pi \int_X \frac{F(u_1, u_2, u_3, u_4)}{u_1^{(b-1)/2} u_2^{(a-1)/2} \sqrt{
            4 u_1 u_2 + a^2 u_1^{b+1} u_2^{a-1} + b^2 u_1^{b-1} u_2^{a+1}
        } } \, d\vol_X.
\end{equation}
Finally, we would like to express $\kappa\, d\vol_X$ as a differential form on $X$ in the $u_k$ coordinates, independent of the parametrization.
To this end, we use the expressions in \eqref{eq:S1-ab_RInvarsInW} to compute
\begin{align*}
    du_1\wedge du_2 \wedge du_3
    &=      -4 a \sqrt{u_1} u_4 \, dw_1\wedge dw_2 \wedge dw_3,
    \\
    du_1\wedge du_2 \wedge du_4
    &=      4 a \sqrt{u_1} u_3 \, dw_1\wedge dw_2 \wedge dw_3,
    \\
    du_1\wedge du_3 \wedge du_4
    &=      2 a^2 u_1^{b + 1/2} u_2^{a-1} \, dw_1\wedge dw_2 \wedge dw_3,
    \\
    du_2\wedge du_3 \wedge du_4
    &=      -2 a b u_1^{b - 1/2} u_2^a \, dw_1\wedge dw_2 \wedge dw_3,
\end{align*}
from which it follows that on $X$, $dw_1\wedge dw_2 \wedge dw_3$ can be expressed as
\[
    \frac{- du_1\wedge du_2 \wedge du_3}
        {4 a \sqrt{u_1} u_4}
    =
    \frac{du_1\wedge du_2 \wedge du_4}
        {4 a \sqrt{u_1} u_3}
    =
    \frac{du_1\wedge du_3 \wedge du_4}
        {2 a^2 u_1^{b + 1/2} u_2^{a-1}}
    =
    \frac{- du_2\wedge du_3 \wedge du_4}
        {2 a b u_1^{b - 1/2} u_2^a}.
\]
We therefore have that \eqref{eq:S1-ab_Int1} is equal to
\[
    2a\pi \int_W g(w_1, w_2, w_3) w_1 \, dw_1 \wedge dw_2 \wedge dw_3
    =
    \int_X F(u_1, u_2, u_3, u_4) \Lambda_X,
\]
where $\Lambda_X$ is given in Equation~\eqref{eq:Main_S1-ab}.

These four expressions for $\Lambda_X$ are of course each distinct as forms on $\R^4$, but any two coincide on the subset of $X$ where both are defined.
The first two are defined on an open dense subset of the principal stratum while the last two are defined on the principal stratum with a singularity at the origin.
Hence, we have that the integrals in Equation~\eqref{eq:S1-ab_Int0}, \eqref{eq:S1-ab_Int1}, and \eqref{eq:S1-ab_Int2} can be expressed as an integral on $X$ as
\[
    \int_X F(u_1, u_2, u_3, u_4) \Lambda_X.
\]
This completes the proof of Theorem~\ref{mainthm:S1-ab}.


\subsection{The symplectic quotient}
\label{subsec:S1-ab_Symp}

Let us now consider the symplectic quotient $Z/\Sp^1$ where $Z = J^{-1}(0)$ is the \emph{shell}, the zero fiber of the moment map.
In this case, the moment map $J\colon V\to \mathfrak{t}^\ast \simeq \R$ is given by
\[
    J(z_1, z_2) = \frac{a}{2} z_1\overline{z_1} - \frac{b}{2} z_2 \overline{z_2}
\]
so that the shell is the set
\begin{equation}
\label{eq:S1-ab_SympShell}
    Z   =   \{ (z_1, z_2) \in V : a z_1\overline{z_1} = b z_2\overline{z_2} \}
        =   \{ (\mathbf{x}, \mathbf{y}) \in V : a (x_1^2 + y_1^2) = b (x_2^2 + y_2^2) \}.
\end{equation}


\subsubsection{The image of the Hilbert embedding and volume form on the shell}
\label{subsubsec:S1-ab_SympImageHilbDS}

In terms of the real invariants, $Z$ is defined by $a u_1 = b u_2$ so that we may use the Hilbert basis $\{ u_1, u_3, u_4 \}$ from \eqref{eq:S1-ab_RInvars},
with $u_2 = a u_1/b$ on the shell, and Hilbert embedding $H_Z\colon Z\to\R^3$ given by $H_Z(x_1,y_1,x_2,y_2) = (u_1, u_3, u_4)$.
The relation \eqref{eq:S1-ab_RRelat} becomes
\begin{equation}
\label{eq:S1-ab_SympRRelat}
    \frac{a^a}{b^a} u_1^{a + b} - u_3^2 - u_4^2 = 0.
\end{equation}
Any point $(u_1, u_3, u_4)\in\R^3$ with $u_1\geq 0$ that satisfies \eqref{eq:S1-ab_SympRRelat} corresponds to a point $(u_1, a u_1/b, u_3, u_4) \in\R^4$
that satisfies \eqref{eq:S1-ab_RRelat}, so the image of the Hilbert embedding $H_Z$ is given by
\begin{equation}
\label{eq:S1-ab_SympImHilb}
    X_Z = \Big\{
    (u_1, u_3, u_4)\in\R^3 :
    \frac{a^a}{b^a} u_1^{a + b} - u_3^2 - u_4^2 = 0,
    \;
    u_1 \geq 0
    \Big\}.
\end{equation}

If $f$ is an $\Sp^1$-invariant smooth function $Z$, by which we mean a smooth $\Sp^1$-invariant function on $V$ restricted to $Z$,
then we consider the integral
\begin{equation}
\label{eq:S1-ab_SympInt0}
    \int_Z f(\mathbf{x}, \mathbf{y}) \, d\vol_Z.
\end{equation}
To compute $d\vol_Z$, we treat $(x_1, y_1, x_2)$ as parameters and use the parameterizations
\[
    (\mathbf{x}, \mathbf{y})
    =
    \Big(x_1, y_1, x_2, \pm\frac{1}{\sqrt{b}} \sqrt{a (x_1^2 + y_1^2) -  bx_2^2}\Big).
\]
The Jacobian is given by
\[
    \frac{\partial (\mathbf{x}, \mathbf{y})}{ \partial (x_1, y_1, x_2)}
    =
    \begin{pmatrix}
           &   I_3   &   
        \\
        \pm\frac{a x_1}{\sqrt{ab (x_1^2 + y_1^2) - b^2 x_2^2}}
            &   \pm\frac{a y_1}{\sqrt{ab (x_1^2 + y_1^2) - b^2 x_2^2}}
            &   \mp\frac{\sqrt{b} x_2}{\sqrt{a(x_1^2 + y_1^2) - b x_2^2}}
    \end{pmatrix},
\]
where $I_3$ is the $3\times 3$ identity matrix. Hence,
\[
    \det \Bigg(\Big(\frac{\partial (\mathbf{x}, \mathbf{y})}{ \partial (x_1, y_1, x_2)}\Big)^T \frac{\partial (\mathbf{x}, \mathbf{y})}{ \partial (x_1, y_1, x_2)}\Bigg)
    =      \frac{ a(a + b) (x_1^2 + y_1^2)}{ab (x_1^2 + y_1^2) - b^2 x_2^2}
\]
does not depend on the sign choice for $y_2$ and
\begin{equation}
\label{eq:S1-ab_SympVolForm}
    d\vol_Z = \sqrt{ \frac{ a(a + b) (x_1^2 + y_1^2)}{ab (x_1^2 + y_1^2) - b^2 x_2^2} } \, dx_1 \wedge dy_1 \wedge dx_2.
\end{equation}


\subsubsection{Reparametrization and the first change of variables}
\label{subsubsec:S1-ab_SympCOV1}

These parametrizations are a natural choice, though they suffer from the fact that they are not $\Sp^1$-invariant; the action of $\Sp^1$ changes the
sign of $y_2$ and hence moves points from one parametrization to the other. This could be remedied using polar coordinates for $V$, but we will rather
reparametrize $Z$ using a fundamental domain analogous to those introduced in Section~\ref{subsec:S1-ab_COV1} but with slightly different coordinates.
Specifically, for each $(\mathbf{x}, \mathbf{y}) \in Z$ with $x_j + i y_j \neq 0$ for $j = 1,2$, we choose the unique orbit representative
$(w_1, 0, w_2, w_3)$ such that $w_1 = \sqrt{x_1^2 + y_1^2} > 0$ and the argument of $w_2 + i w_3$ is contained in $[0, 2\pi/a)$. We parametrize
this set using $w_1$ and $\alpha = \arg (w_2 + i w_3) \in [0, 2\pi/a)$.

It will be convenient to set $R = a/b$.
Noting that $(w_1, 0, w_2, w_3)\in Z$ implies via the description \eqref{eq:S1-ab_SympShell} that
\[
    \abs{w_2 + i w_3}   =   \sqrt{R} w_1,
\]
we express the fundamental domain as
\begin{equation}
\label{eq:S1-ab_SympFundDom}
    \{(w_1, 0, w_2, w_3)\in\R^4 : w_1 \geq 0 \}  = \Big\{ \Big( w_1, 0, \sqrt{R} w_1 \cos\alpha, \sqrt{R} w_1 \sin\alpha \Big) \in\R^4 :
    w_1 \geq 0, \;
    \alpha\in [0, 2\pi/a) \Big\}.
\end{equation}
Choosing the unique $\theta\in(-\pi,\pi]$ such that $e^{i\theta}\ast(w_1, 0, w_2, w_3) = (\mathbf{x}, \mathbf{y})$, we have
\begin{equation}
\label{eq:S1-ab_SympCOV1}
    (\mathbf{x}, \mathbf{y})
    =
    \big(   w_1 \cos(a\theta),
            - w_1 \sin(a\theta),
            \sqrt{R} w_1 \cos(\alpha + b\theta),
            \sqrt{R} w_1 \sin(\alpha + b\theta)
    \big).
\end{equation}
Recalling that $\theta$ is oriented negatively as in Section~\ref{subsec:S1-ab_COV1}, the Jacobian is given by
\[
    \frac{\partial (\mathbf{x}, \mathbf{y})}{ \partial (\theta, w_1, \alpha)}
    =
    \begin{pmatrix}
        - a w_1\sin(a\theta)     &      \cos(a\theta)   &      0
        \\
        - a w_1\cos(a\theta)     &      - \sin(a\theta) &       0
        \\
        -b \sqrt{R} w_1 \sin(\alpha + b\theta)
            &   \sqrt{R}\cos(\alpha + b\theta)
            &   - \sqrt{R} w_1 \sin(\alpha + b\theta)
        \\
        b \sqrt{R} w_1 \cos(\alpha + b\theta)
            &   \sqrt{R}\sin(\alpha + b\theta)
            &   \sqrt{R} w_1 \cos(\alpha + b\theta)
    \end{pmatrix},
\]
and
\[
    \det\Bigg(\Big(\frac{\partial (\mathbf{x}, \mathbf{y})}{ \partial (\theta, w_1,\alpha)}\Big)^T \frac{\partial (\mathbf{x}, \mathbf{y})}{ \partial (\theta,w_1,\alpha)}\Bigg)
    =
    R^2 a (a + b) w_1^4
\]
is independent of $\alpha$ and $\theta$, so that in terms of this parametrization,
\begin{equation}
\label{eq:S1-ab_SympVolFormW}
    d\vol_Z =  R \sqrt{a(a + b)} w_1^2 \, d\theta \wedge dw_1 \wedge d\alpha,
    \qquad
    R = \frac{a}{b}.
\end{equation}
We compose $f$ with the change of variables \eqref{eq:S1-ab_SympCOV1} to express $f(x_1,y_1,x_2,y_2) = g(w_1, \alpha)$ for a smooth function $g$ and
recall that invariance of $f$ implies that $g$ does not depend on $\theta$.
We therefore can express the integral \eqref{eq:S1-ab_SympInt0} as an integral over the fundamental domain in terms of this parametrization, given by
\begin{align}
    \nonumber
    &\int_{\substack{\theta\in(-\pi,\pi] \\ w_1 \geq 0 \\ \alpha\in[0,2\pi/a) }} g(w_1, \alpha)
        R \sqrt{a(a + b)} w_1^2
        \, d\theta \wedge dw_1 \wedge d\alpha
    \\
    \nonumber
    &\qquad=
    R \sqrt{a(a + b)}  \int_{\substack{w_1 \geq 0 \\ \alpha\in[0,2\pi/a)}} g(w_1, \alpha)
        w_1^2
        \, dw_1 \wedge  d\alpha \int_{-\pi}^\pi d\theta
    \\
    \label{eq:S1-ab_SympInt1}
    &\qquad=
    2\pi R \sqrt{a(a + b)}  \int_{\substack{w_1 \geq 0 \\ \alpha\in[0,2\pi/a)}} g(w_1, \alpha)
        w_1^2 \, dw_1 \wedge  d\alpha .
\end{align}


\subsubsection{The second change of variables}
\label{subsubsec:S1-ab_SympCOV2}

We now turn to expressing the integral \eqref{eq:S1-ab_SympInt1} as an integral over the image $X_Z$ of the Hilbert embedding $H_Z$.
For a point in the fundamental domain \eqref{eq:S1-ab_SympFundDom}, we have
\begin{equation}
\label{eq:S1-ab_SympRInvarsInW}
\begin{split}
    u_1 &=      w_1^2
    \qquad\qquad
    u_3 =       R^{a/2} w_1^{a + b} \cos(\alpha),
    \\
    u_4 &=      R^{a/2} w_1^{a + b} \sin(\alpha).
\end{split}
\end{equation}
To compute the volume form on $X_Z$, we have
\[
    \frac{\partial (u_1, u_3, u_4)}{ \partial (w_1, \alpha)}
    =
    \begin{pmatrix}
        2w_1    &   0
        \\
        (a+b) w_1^{a + b - 1} R^{a/2} \cos(\alpha)
            &   - w_1^{a + b} R^{a/2} \sin(\alpha)
        \\
        (a+b) w_1^{a + b - 1} R^{a/2} \sin(\alpha)
            &   w_1^{a + b} R^{a/2} \cos(\alpha)
    \end{pmatrix},
\]
so that
\begin{align*}
    d\vol_{X_Z}
    &=      \sqrt{\det \Bigg(\Big(\frac{\partial (u_1, u_3, u_4)}{ \partial (w_1, \alpha)}\Big)^T \frac{\partial (u_1, u_3, u_4)}{ \partial (w_1, \alpha)}\Bigg)}
                \, dw_1 \wedge d\alpha
    \\&=    R^{a/2} u_1^{(a + b + 1)/2}  \sqrt{ 4 + R^a (a + b)^2 u_1^{a + b - 2} } \, dw_1 \wedge d\alpha.
\end{align*}
Letting $F$ be a smooth function on $\R^3$ with $F(u_1, u_3, u_4) = g(w_1, \alpha) = f(\mathbf{x}, \mathbf{y})$, given by applying \cite[Thrm.~1]{GWSdiff} to the
smooth function on $V$ whose restriction to $Z$ is $f$, the equation
\[
     2\pi R \sqrt{a(a + b)}  \int_{\substack{w_1 \geq 0 \\ \alpha\in[0,2\pi/a)}} g(w_1, \alpha)
        w_1^2 \, dw_1 \wedge  d\alpha
     =
     \int_{X_Z} F(u_1, u_3, u_4) \kappa \, d\vol_{X_Z}
\]
implies
\[
    \kappa
    =
    \frac{ 2\pi \sqrt{a(a + b)} }
    { R^{(a-2)/2} u_1^{(a + b - 1)/2} \sqrt{ 4 + R^a (a + b)^2 u_1^{a + b - 2} } }
\]
so we can write \eqref{eq:S1-ab_SympInt1} as
\begin{equation}
\label{eq:S1-ab_SympInt2}
    2\pi \sqrt{\frac{a(a+b)}{R^{a-2}}} \int_{X_Z}
        \frac{ F(u_1, u_3, u_4) }{ u_1^{(a + b - 1)/2} \sqrt{4 + R^a(a + b)^2 u_1^{a + b - 2} } }
    \, d\vol_{X_Z}.
\end{equation}
Finally, we compute from \eqref{eq:S1-ab_SympRInvarsInW} that
\begin{align*}
    du_1\wedge du_3
    &=      -2 R^{a/2} w_1^{a+b+1}\sin(\alpha) \, dw_1 \wedge d\alpha
    &&=     -2 \sqrt{u_1} u_4 \, dw_1 \wedge d\alpha
    \\
    du_1\wedge du_4
    &=      2 R^{a/2} w_1^{a+b+1}\cos(\alpha) \, dw_1 \wedge d\alpha
    &&=     2 \sqrt{u_1} u_3 \, dw_1 \wedge d\alpha
    \\
    du_3\wedge du_4
    &=      (a + b) R^a w_1^{2a + 2b - 1}  \, dw_1 \wedge d\alpha
    &&=     (a + b) R^a u_1^{a + b - 1/2}  \, dw_1 \wedge d\alpha
\end{align*}
so that
\[
    dw_1 \wedge d\alpha
    =
    -\frac{du_1\wedge du_3}{ 2 \sqrt{u_1} u_4 }
    =
    \frac{du_1\wedge du_4}{ 2 \sqrt{u_1} u_3 }
    =
    \frac{du_3\wedge du_4}{ (a + b) R^a u_1^{a + b - 1/2} },
\]
and the integrals in \eqref{eq:S1-ab_SympInt0} and \eqref{eq:S1-ab_SympInt1} are given by
\[
    \int_{X_Z} F(u_1, u_3, u_4) \Lambda_{X_Z}
\]
where $\Lambda_{X_Z}$ is given in Equation~\eqref{eq:Main_S1-ab_Symp}.
Once again, the first two expressions for $\Lambda_{X_Z}$ are defined on an open dense subset of the principal stratum while
the last is defined throughout the principal stratum with a singularity at the point fixed by $\Sp^1$.
This completes the proof of Theorem~\ref{mainthm:S1-ab_Symp}.


\section{Circle representations with weight vector $(-1, a, a)$}
\label{sec:S1-1aa}

We now consider the case where $V\simeq\C^3$ is the representation of $G = \Sp^1$ with weight vector $(-1, a, a)$ where $a\in\Z$ is positive;
see the beginning of Section~\ref{sec:S1-ab}.
Following the outline of Section~\ref{sec:S1-ab}, we describe the Hilbert embedding and invariants of the representation $V$ in
Section~\ref{subsec:S1-1aa_ImageHilb}, compute the Hilbert measure for $V$ in Sections~\ref{subsec:S1-1aa_COV1} and \ref{subsec:S1-1aa_COV2}
and hence prove Theorem~\ref{mainthm:S1-1aa}, and compute the Hilbert measure for the corresponding symplectic quotient at level zero in
Section~\ref{subsec:S1-1aa_Symp}, proving Theorem~\ref{mainthm:S1-1aa_Symp}.


\subsection{The image of the Hilbert embedding}
\label{subsec:S1-1aa_ImageHilb}

In complex coordinates $(z_1, z_2, z_3)$ for $V$, a minimal set of generating invariants consists of
\begin{equation}
\label{eq:S1-1aa_CInvars}
\begin{aligned}
    m_1 &=      z_1 \overline{z_1},
    \qquad&
    m_2 &=      z_2 \overline{z_2},
    \qquad&
    m_3 &=      z_3 \overline{z_3},
    \\
    p_1 &=      z_1^a z_2,
    &
    p_2 &=      z_1^a z_3,
    &
    p_3 &=      z_2 \overline{z_3},
    \\
    q_1 &=      \overline{z_1}^a \overline{z_2},
    &
    q_2 &=      \overline{z_1}^a \overline{z_3},
    &
    q_3 &=      z_3 \overline{z_2},
\end{aligned}
\end{equation}
and the ideal of relations is generated by the nine $2\times 2$ minors of the matrix
\begin{equation}
\label{eq:S1-1aa_CRelat}
    \begin{pmatrix}
        m_1^a   &   p_1 &   p_2
        \\
        q_1     &   m_2 &   q_3
        \\
        q_2     &   p_3 &   m_3
    \end{pmatrix}
    =
    \begin{pmatrix}
        z_1^a
        \\
        \overline{z_2}
        \\
        \overline{z_3}
    \end{pmatrix}
    \begin{pmatrix}
        \overline{z_1}^a
        &
        z_2
        &
        z_3
    \end{pmatrix}.
\end{equation}
Letting $z_j = x_j + i y_j$ for $j = 1,2,3$ and $(\mathbf{x},\mathbf{y}) = (x_1, y_1, x_2, y_2, x_3, y_3)$ be real coordinates for $V$,
the invariants $m_j$ are real and given by $m_j = x_j^2 + y_j^2$.
The remaining real invariants are given by the real and imaginary parts of the $p_j$, which again coincide up to sign with the real and imaginary parts of the $q_j$.
Hence the real Hilbert basis is given by
\begin{equation}
\label{eq:S1-1aa_RInvars}
\begin{split}
    u_1 &=      x_1^2 + y_1^2,
    \qquad
    u_2 =       x_2^2 + y_2^2,
    \qquad
    u_3 =       x_3^2 + y_3^2,
    \\
    u_4 &=      \Real p_1
        =       \frac{1}{2}\big( (x_1 + i y_1)^a + (x_1 - i y_1)^a \big)x_2 - \frac{1}{2i}\big( (x_1 + i y_1)^a - (x_1 - i y_1)^a \big)y_2
    \\
    u_5 &=      \Imaginary p_1
        =       \frac{1}{2}\big( (x_1 + i y_1)^a + (x_1 - i y_1)^a \big)y_2 + \frac{1}{2i}\big( (x_1 + i y_1)^a - (x_1 - i y_1)^a \big)x_2
    \\
    u_6 &=      \Real p_2
        =       \frac{1}{2}\big( (x_1 + i y_1)^a + (x_1 - i y_1)^a \big)x_3 - \frac{1}{2i}\big( (x_1 + i y_1)^a - (x_1 - i y_1)^a \big)y_3
    \\
    u_7 &=      \Imaginary p_2
        =       \frac{1}{2}\big( (x_1 + i y_1)^a + (x_1 - i y_1)^a \big)y_3 + \frac{1}{2i}\big( (x_1 + i y_1)^a - (x_1 - i y_1)^a \big)x_3
    \\
    u_8 &=      \Real p_3
        =       x_2 x_3 + y_2 y_3
    \qquad
    u_9 =      \Imaginary p_3
        =       x_3 y_2 - x_2 y_3.
\end{split}
\end{equation}
It follows that the Hilbert embedding $H\colon V\to\R^9$ is given by
$H(\mathbf{x}, \mathbf{y}) = (u_1, \ldots, u_9)$. From Equation~\eqref{eq:S1-1aa_CRelat}, one computes that the real relations are generated by
\begin{equation}
\label{eq:S1-1aa_RRelat}
\begin{split}
    u_2 u_3 - u_8^2 - u_9^2, \quad
    u_2 u_6 - u_4 u_8 - u_5 u_9, \quad
    u_2 u_7 - u_5 u_8 + u_4 u_9,
    \\
    u_3 u_4 - u_6 u_8 + u_7 u_9, \quad
    u_3 u_5 - u_7 u_8 - u_6 u_9, \quad
    u_1^a u_2 - u_4^2 - u_5^2,
    \\
    u_4 u_6 + u_5 u_7 - u_1^a u_8, \quad
    u_5 u_6 - u_4 u_7 - u_1^a u_9, \quad
    u_1^a u_3 - u_6^2 - u_7^2.
\end{split}
\end{equation}

We claim the following.

\begin{lemma}
\label{lem:S1-1aa_ImHilb}
The image $X$ of the Hilbert embedding associated to the invariants \eqref{eq:S1-1aa_RInvars} is the set of $(u_1,\ldots,u_9)\in\R^9$ satisfying the relations
\eqref{eq:S1-1aa_RRelat} and the inequalities $u_1 \geq 0$, $u_2 \geq 0$, and $u_3\geq 0$.
\end{lemma}
\begin{proof}
From the definition of the $u_j$, it is clear that any $(u_1,\ldots,u_9)\in X$ satisfies the three given inequalities as well as the relations \eqref{eq:S1-1aa_RRelat}.
So assume $(u_1,\ldots,u_9)\in\R^9$ satisfies the given inequalities and define the complex numbers
\[
    P_1 =   u_4 + i u_5,
    \quad
    P_2 =   u_6 + i u_7,
    \quad
    P_3 =   u_8 + i u_9.
\]
Applying the sixth, ninth, first, seventh, and eighth relations, respectively, we have
\begin{align}
    \label{eq:S1-1aa_ImHilbZ2}
    \abs{P_1}^2
    &=  u_4^2 + u_5^2
    =   u_1^a u_2,
    \\
    \label{eq:S1-1aa_ImHilbZ3}
    \abs{P_2}^2
    &=  u_6^2 + u_7^2
    =   u_1^a u_3,
    \\
    \label{eq:S1-1aa_ImHilbZ4}
    \abs{P_3}^2
    &=  u_8^2 + u_9^2
    =   u_2 u_3,
    \\
    \label{eq:S1-1aa_ImHilbZ2Z3bar}
    P_1\overline{P_2}
    &=  u_4 u_6 + u_5 u_7 + i (u_5 u_6 - u_4 u_7)
    =   u_1^a u_8 + i u_1^a u_9
    =   u_1^a P_3.
\end{align}
If $u_1 = 0$, then \eqref{eq:S1-1aa_ImHilbZ2} and \eqref{eq:S1-1aa_ImHilbZ3} imply that $P_1 = P_2 = 0$.
For any $z_2, z_3 \in \C$ such that $\abs{z_2} = \sqrt{u_2}$ and $\abs{z_3} = \sqrt{u_3}$,
we have by \eqref{eq:S1-1aa_ImHilbZ4} that $\abs{z_2 z_3} = \abs{z_2}\abs{z_3} = \abs{P_3}$, so we may choose such $z_2, z_3$ that satisfy $z_2\overline{z_3} = P_3$.
Then $(0,z_2,z_3)$ are the complex coordinates of a point $(x_1,y_1,x_2,y_2,x_3,y_3)\in V$ whose image under $H$ is $(u_1,\ldots,u_9)$.

Now assume $u_1\neq 0$, choose $z_1\in\C$ such that $\abs{z_1}^2 = u_1$, and define
\[
    z_2 = \frac{P_1}{z_1^a},
    \quad
    z_3 = \frac{P_2}{z_1^a}.
\]
Then by \eqref{eq:S1-1aa_ImHilbZ2}, $\abs{z_2}^2 = u_2$, by \eqref{eq:S1-1aa_ImHilbZ3}, $\abs{z_3}^2 = u_3$, and by
\eqref{eq:S1-1aa_ImHilbZ2Z3bar}, $z_2\overline{z_3} = P_3$.
Then $(z_1,z_2,z_3)$ are the complex coordinates of a point in $V$ whose image under $H$ is $(u_1,\ldots,u_9)$, completing the proof.
\end{proof}


\subsection{The first change of variables}
\label{subsec:S1-1aa_COV1}

Let $f(\mathbf{x}, \mathbf{y})$ be a smooth, real-valued, $\Sp^1$-invariant integrable function on $V$. We transform the integral
\begin{equation}
\label{eq:S1-1aa_Int0}
    \int_V f(\mathbf{x}, \mathbf{y}) \, dx_1 \wedge dy_1 \wedge dx_2\wedge dy_2 \wedge dx_3\wedge dy_3
\end{equation}
using the same approach as in Section~\ref{subsec:S1-ab_COV1}. A fundamental domain in $V$ is given by points of the form $(w_1, 0, w_2, w_3, w_4, w_5) \in V$
with $w_1\geq 0$, and we use this domain and the $\Sp^1$-action to parametrize $V$ as
\begin{multline}
\label{eq:S1-1aa_COV1}
    (\mathbf{x}, \mathbf{y})
    =
    \big(   w_1 \cos\theta,
            - w_1 \sin\theta,
            w_2 \cos(a\theta) - w_3 \sin(a\theta),
            w_2 \sin(a\theta) + w_3 \cos(a\theta),
    \\
            w_4 \cos(a\theta) - w_5 \sin(a\theta),
            w_4 \sin(a\theta) + w_5 \cos(a\theta)
    \big).
\end{multline}
The Jacobian $\frac{\partial (\mathbf{x}, \mathbf{y})}{ \partial (\theta, w_1, w_2, w_3, w_4, w_5)}$ has a block form similar to that of
\eqref{eq:S1-ab_COV1Jacobian} with diagonal blocks
\[
        \begin{pmatrix}
        - w_1\sin\theta   &   \cos\theta
        \\
        - w_1\cos\theta   &   - \sin\theta
        \end{pmatrix},
        \qquad
        \begin{pmatrix}
        \cos(a \theta)   &  - \sin(a\theta)
        \\
        \sin(a \theta)  &   \cos(a\theta)
        \end{pmatrix},
        \qquad
        \begin{pmatrix}
        \cos(a \theta)   &  - \sin(a\theta)
        \\
        \sin(a \theta)  &   \cos(a\theta)
        \end{pmatrix},
\]
hence a determinant of $w_1$. Then
\[
    dx_1\wedge dy_1 \wedge dx_2 \wedge dy_2 \wedge dx_3 \wedge dy_3
    =
    w_1 d\theta \wedge dw_1 \wedge dw_2 \wedge dw_3 \wedge dw_4 \wedge dw_5.
\]
Letting $g(w_1,\ldots,w_5)$ denote the composition of $f$ with \eqref{eq:S1-1aa_COV1}, which is independent of $\theta$ by $\Sp^1$-invariance,
and denoting
\[
    d\mathbf{w} =   dw_1 \wedge \cdots \wedge dw_5,
\]
we have that the integral \eqref{eq:S1-1aa_Int0} is equal to
\begin{align}
    \nonumber
    &\int_{(-\pi,\pi]\times W} g(w_1, \ldots, w_5) w_1 \, d\theta \wedge d\mathbf{w}
    \\
    \nonumber
    &\qquad=    \int_{-\pi}^{\pi} d\theta \int_W g(w_1, \ldots, w_5) w_1 \, d\mathbf{w}
    \\
    \label{eq:S1-1aa_Int1}
    &\qquad=    2\pi \int_W g(w_1, \ldots, w_5) w_1 \, d\mathbf{w},
\end{align}
where $W$ is the fundamental domain
\[
    W = \{ (w_1, 0, w_2, w_3, w_4, w_5) \in V : w_1 > 0 \}.
\]


\subsection{The second change of variables}
\label{subsec:S1-1aa_COV2}

On the fundamental domain $W$, the real invariants in \eqref{eq:S1-1aa_RInvars} simplify to
\begin{equation}
\label{eq:S1-1aa_RInvarsInW}
\begin{aligned}
    u_1 &=      w_1^2,
    &
    u_2 &=      w_2^2 + w_3^2,
    &u_3 &=     w_4^2 + w_5^2,
    \\
    u_4 &=      w_1^a w_2,
    \qquad
    &
    u_5 &=      w_1^a w_3,
    \qquad
    &
    u_6 &=      w_1^a w_4,
    \\
    u_7 &=      w_1^a w_5,
    &
    u_8 &=      w_2 w_4 + w_3 w_5,
    &
    u_9 &=      w_3 w_4 - w_2 w_5.
\end{aligned}
\end{equation}
We compute the volume form $d\vol_X$ on $X$ by first computing the Jacobian
\[
    \frac{\partial\mathbf{u}}{\partial\mathbf{w}} = \frac{\partial (u_1,\ldots,u_9)}{\partial (w_1,\ldots,w_5)}
    =
    \begin{pmatrix}
        2 w_1 & 0 & 0 & 0 & 0 \\
        0 & 2 w_2 & 2 w_3 & 0 & 0 \\
        0 & 0 & 0 & 2 w_4 & 2 w_5 \\
        a w_1^{a-1} w_2  & w_1^a & 0 & 0 & 0 \\
        a w_1^{a-1} w_3 & 0 & w_1^a & 0 & 0 \\
        a w_1^{a-1} w_4 & 0 & 0 & w_1^a & 0 \\
        a w_1^{a-1} w_5 & 0 & 0 & 0 & w_1^a \\
        0 & w_4 & w_5 & w_2 & w_3 \\
        0 & -w_5 & w_4 & w_3 & -w_2
    \end{pmatrix},
\]
and then
\[
    \det \Bigg( \Big(\frac{\partial\mathbf{u}}{\partial\mathbf{w}}\Big)^T \frac{\partial\mathbf{u}}{\partial\mathbf{w}} \Bigg) 
    =
    4 u_1^a (u_1^a + u_2 + u_3) L
\]
where
\[
    L
    =   \Big(
        4 u_1 (u_2^2 + 4 u_2 u_3 + u_3^2) + u_1^{2a - 1} \big( u_1^2 + a^2 (u_2^2 + u_2 u_3 + u_3^2)\big)
        + u_1^{a - 1} (u_2 + u_3) \big(5 u_1^2 + a^2 (u_2^2 + 4 u_2 u_3 + u_3^2)\big)
        \Big).
\]
Hence
\[
    d\vol_X
    =
    2 u_1^{a/2} \sqrt{ (u_1^a + u_2 + u_3) L }
        \, d\mathbf{w}.
\]
For a smooth function $F$ with $F(u_1, \ldots, u_9) = g(w_1, \ldots, w_5) = f(\mathbf{x}, \mathbf{y})$, whose existence is guaranteed by
\cite[Thrm.~1]{GWSdiff}, we solve
\[
    2\pi \int_W g(w_1, \ldots, w_5) w_1 \, d\mathbf{w}
    =
    \int_X F(u_1, \ldots, u_9) \kappa \, d\vol_X
\]
to obtain
\[
    \kappa
    =
    \frac{\pi}{
        u_1^{(a-1)/2} \sqrt{ (u_1^a + u_2 + u_3) L }
        },
\]
and then \eqref{eq:S1-1aa_Int1} can be written in terms of the parametrization \eqref{eq:S1-1aa_RInvarsInW} as
\begin{equation}
\label{eq:S1-1aa_Int2}
    \pi \int_X \frac{F(u_1, \ldots, u_9)}
                    {u_1^{(a-2)/2} \sqrt{ (u_1^a + u_2 + u_3) L }} \, d\vol_X .
\end{equation}

Unlike the situations in Sections~\ref{subsec:S1-ab_COV2} and \ref{subsubsec:S1-ab_SympCOV2}, it is infeasible to give a
complete list of expressions of the Hilbert measure based on different coordinate choices. To give two expressions,
we compute from \eqref{eq:S1-1aa_RInvarsInW} that
\begin{align*}
    du_1\wedge du_2 \wedge du_3 \wedge du_4 \wedge du_7
    =      8 w_1 u_5 u_6 \, d\mathbf{w}
    \qquad\text{ and }\qquad
    du_1\wedge du_2 \wedge du_3 \wedge du_5 \wedge du_6
    =      8 w_1 u_4 u_7 \, d\mathbf{w}.
\end{align*}
It follows that $\Lambda_X$ is as given in Equation~\eqref{eq:Main_S1-1aa}, with each expression defined on a dense subset of the principal stratum,
completing the proof of Theorem~\ref{mainthm:S1-1aa}.


\subsection{The symplectic quotient}
\label{subsec:S1-1aa_Symp}

Turning to the symplectic quotient, the moment map in this case is given by
\[
    J(z_1, z_2, z_3) = \frac{1}{2}z_1\overline{z_1} - \frac{a}{2} z_2 \overline{z_2}  - \frac{a}{2} z_3 \overline{z_3},
\]
hence the shell is
\begin{equation}
\label{eq:S1-1aa_SympShell}
    Z   =   \{ (z_1, z_2, z_3) \in V : z_1\overline{z_1} = a (z_2\overline{z_2} + z_3\overline{z_3}) \}
        =   \{ (\mathbf{x}, \mathbf{y}) \in V : x_1^2 + y_1^2 = a (x_2^2 + y_2^2 + x_3^2 + y_3^2) \}.
\end{equation}

As $Z$ is defined by $u_1 = a(u_2 + u_3)$, we use the Hilbert basis $\{ u_1, u_2, u_4, \ldots, u_9 \}$ from \eqref{eq:S1-1aa_RInvars}
with $u_3 = u_1/a - u_2$ on the shell, and Hilbert embedding $H_Z\colon Z\to\R^8$ given by $H_Z(\mathbf{x}, \mathbf{y}) = (u_1, u_2, u_4, \ldots, u_9)$.
As in Section~\ref{subsubsec:S1-ab_SympImageHilbDS}, the image $X_Z$ of $H_Z$ is given by the set of points $(u_1, u_2, u_4, \ldots, u_9)\in\R^8$
such that $(u_1, u_2, u_1/a - u_2, u_4 \ldots,u_9)\in X$, i.e., satisfying \eqref{eq:S1-1aa_RRelat} with $u_3 = u_1/a - u_2$ as well as $u_1, u_3 \geq 0$.

We are once again interested in the integral
\begin{equation}
\label{eq:S1-1aa_SympInt0}
    \int_Z f(\mathbf{x},\mathbf{y}) \, d\vol_Z
\end{equation}
for an $\Sp^1$-invariant smooth function $Z$. We can parametrize $Z$ as
\[
    (\mathbf{x},\mathbf{y})
    =
    \Big(x_1, y_1, x_2, y_2, x_3, \pm\frac{1}{\sqrt{a}} \sqrt{x_1^2 + y_1^2 -  a(x_2^2 + y_2^2 + x_3^2)}\Big),
\]
with Jacobian
\[
    \frac{\partial (\mathbf{x}, \mathbf{y})}{ \partial (x_1, y_1, x_2, y_2, x_3)}
    =
    \begin{pmatrix}
           &&   I_5   &&
        \\
        \pm \frac{x_1}{\sqrt{a} B} & \pm \frac{y_1}{\sqrt{a} B} & \mp \frac{\sqrt{a} x_2}{B} & \mp \frac{\sqrt{a} y_2}{B} & \mp \frac{\sqrt{a} x_3}{B}
    \end{pmatrix},
    \qquad
    B = \sqrt{x_1^2 + y_1^2 -  a(x_2^2 + y_2^2 + x_3^2)},
\]
where $I_5$ is the identity matrix, and then
\[
    \det \Bigg(\Big(\frac{\partial (\mathbf{x}, \mathbf{y})}{ \partial (x_1, y_1, x_2, y_2, x_3)}\Big)^T
        \frac{\partial (\mathbf{x}, \mathbf{y})}{ \partial (x_1, y_1, x_2, y_2, x_3)}\Bigg)
    =      \frac{ (1 + a) (x_1^2 + y_1^2)}{a (x_1^2 + y_1^2) - a^2 (x_2^2 + y_2^2 + x_3^2)}
\]
is nonnegative on $Z$, does not depend on the sign choice, and
\begin{equation}
\label{eq:S1-1aa_SympVolForm}
    d\vol_Z = \sqrt{ \frac{ (1 + a) (x_1^2 + y_1^2)}{a (x_1^2 + y_1^2) - a^2 (x_2^2 + y_2^2 + x_3^2)} } \, dx_1 \wedge dy_1 \wedge dx_2 \wedge dy_2 \wedge dx_3.
\end{equation}

As in Section~\ref{subsubsec:S1-ab_SympCOV1}, we parametrize $Z$ away from the origin by choosing from each orbit the unique point of the form
$(w_1, 0, w_2, w_3, w_4, w_5)$ with $w_1 \geq 0$. We set $\alpha = \arg (w_4 + i w_5) \in (-\pi, \pi]$, note that \eqref{eq:S1-1aa_SympShell} implies
\[
    w_4^2 + w_5^2   =   \frac{w_1^2}{a} - w_2^2 - w_3^2,
\]
and express points in the fundamental domain as
\begin{equation}
\label{eq:S1-1aa_SympFundDom}
    (w_1, 0, w_2, w_3, w_4, w_5) = \Big( w_1, 0, w_2, w_3, \sqrt{w_1^2/a - w_2^2 - w_3^2}\cos\alpha, \sqrt{w_1^2/a - w_2^2 - w_3^2}\sin\alpha \Big),
\end{equation}
where $w_1 \geq 0$ and $\alpha\in (-\pi, \pi]$.
Choosing the unique $\theta\in(-\pi,\pi]$ such that $e^{i\theta}\ast(w_1, 0, w_2, w_3, w_4, w_5) = (\mathbf{x}, \mathbf{y})$, we parametrize $Z$
using parameters $(\theta,w_1,w_2,w_3,\alpha)$ as
\begin{multline}
\label{eq:S1-1aa_SympCOV1}
    (\mathbf{x}, \mathbf{y})
    =
    \big(   w_1 \cos\theta,
            - w_1 \sin\theta,
            w_2 \cos(a\theta) - w_3 \sin(a\theta),
            w_2 \sin(a\theta) + w_3 \cos(a\theta), \\
            \sqrt{w_1^2/a - w_2^2 - w_3^2}\cos(\alpha + a\theta),
            \sqrt{w_1^2/a - w_2^2 - w_3^2}\sin(\alpha + a\theta),
    \big).
\end{multline}
Recalling the orientation of $\theta$ from Section~\ref{subsec:S1-1aa_COV1}, and abbreviating
$\mathbf{w}_{123} = (w_1, w_2, w_3)$ and $d\mathbf{w}_{123} = dw_1 \wedge dw_2 \wedge dw_3$, we compute
\[
    \det\Bigg(\Big(\frac{\partial (\mathbf{x}, \mathbf{y})}{ \partial (\theta, \mathbf{w}_{123}, \alpha)}\Big)^T
                \frac{\partial (\mathbf{x}, \mathbf{y})}{ \partial (\theta, \mathbf{w}_{123}, \alpha)}\Bigg)
    =
    \frac{1 + a}{a^4} w_1^4,
\]
hence
\begin{equation}
\label{eq:S1-1aa_SympVolFormW}
    d\vol_Z =  \frac{\sqrt{1 + a}}{a^2} w_1^2
                \, d\theta \wedge d\mathbf{w}_{123} \wedge d\alpha.
\end{equation}
Then, composing $f$ with \eqref{eq:S1-1aa_SympCOV1} to obtain a smooth $g$ such that $f(\mathbf{x},\mathbf{y}) = g(w_1,w_2,w_3,\alpha)$,
and integrating over $\theta$ as in Section~\ref{subsec:S1-ab_COV1}, \ref{subsubsec:S1-ab_SympCOV1}, and \ref{subsec:S1-1aa_COV1},
we express \eqref{eq:S1-1aa_SympInt0} as
\begin{equation}
    \label{eq:S1-1aa_SympInt1}
    \frac{2\pi\sqrt{1 + a}}{a^2}  \int\limits_{\substack{w_1 \geq 0 \\ w_2, w_3\in\R \\ \alpha\in(-\pi,\pi)}} g(\mathbf{w}_{123},\alpha)
        w_1^2 \, d\mathbf{w}_{123} \wedge d\alpha.
\end{equation}

For a point in the fundamental domain \eqref{eq:S1-1aa_SympFundDom} corresponding to parameters $(\mathbf{w}_{123}, \alpha)$, the real invariants
\eqref{eq:S1-1aa_RInvarsInW} are given in terms of the parametrization \eqref{eq:S1-1aa_SympCOV1} by
\begin{equation}
\label{eq:S1-1aa_SympRInvarsInW}
\begin{aligned}
    u_1 &=      w_1^2,
    \qquad\qquad
    u_2 =       w_2^2 + w_3^2,
    \qquad\qquad
    u_4 =       w_1^a w_2,
    \qquad\qquad
    u_5 =       w_1^a w_3,
    \\
    u_6 &=      w_1^a \sqrt{w_1^2/a - w_2^2 - w_3^2}\cos\alpha,
    \qquad\qquad
    u_7 =       w_1^a \sqrt{w_1^2/a - w_2^2 - w_3^2}\sin\alpha,
    \\
    u_8 &=      w_2 \sqrt{w_1^2/a - w_2^2 - w_3^2}\cos\alpha + w_3 \sqrt{w_1^2/a - w_2^2 - w_3^2}\sin\alpha,
    \\
    u_9 &=      w_3 \sqrt{w_1^2/a - w_2^2 - w_3^2}\cos\alpha - w_2\sqrt{w_1^2/a - w_2^2 - w_3^2}\sin\alpha.
\end{aligned}
\end{equation}
From these, we compute
\begin{align*}
    d\vol_{X_Z}
    &=
    \sqrt{\det\Bigg(\Big(\frac{\partial (u_1, u_2, u_4,\ldots,u_9)}{ \partial (\mathbf{w}_{123}, \alpha)}\Big)^T
        \frac{\partial (u_1, u_2, u_4,\ldots,u_9)}{ \partial (\mathbf{w}_{123}, \alpha)}\Bigg) }
    \, d\mathbf{w}_{123}\wedge d\alpha
    \\&=
    \frac{u_1^{(a + 1)/2} }{a^2} \sqrt{ (u_1 + a u_1^a)\big(a(a + 1)^2 u_1^{2a} + (5a^2 + 2a + 1) u_1^{a + 1} + 4a u_1^a u_2 + 4a u_1^2 + 4a u_2^2 \big) }
     \, d\mathbf{w}_{123}\wedge d\alpha.
\end{align*}
Applying \cite[Thrm.~1]{GWSdiff}, for a smooth function $F$ on $\R^8$ with $F(u_1, u_2, u_4,\ldots, u_9) = g(\mathbf{w}_{123}, \alpha) = f(\mathbf{x}, \mathbf{y})$
on $Z$,
\begin{align*}
    &\frac{2\pi\sqrt{1 + a}}{a^2}  \int\limits_{\substack{w_1 \geq 0 \\ w_2, w_3\in\R \\ \alpha\in(-\pi,\pi)}} g(\mathbf{w}_{123},\alpha)
        w_1^2 \, d\mathbf{w}_{123} \wedge d\alpha
    \\&\qquad=
    \frac{2\pi\sqrt{1 + a}}{a^2}  \int\limits_{\substack{w_1 \geq 0 \\ w_2, w_3\in\R \\ \alpha\in(-\pi,\pi)}} g(\mathbf{w}_{123},\alpha)
            u_1 \, d\mathbf{w}_{123} \wedge d\alpha
    \\&\qquad=
    \int_{X_Z} F(u_1, u_2, u_4,\ldots, u_9) \kappa \, d\vol_{X_Z}
\end{align*}
where
\[
    \kappa
    =
    \frac{2\pi\sqrt{1 + a}}
        {u_1^{(a - 1)/2}
            \sqrt{ (u_1 + a u_1^a)\big(a(a + 1)^2 u_1^{2a} + (5a^2 + 2a + 1) u_1^{a + 1} + 4a u_1^a u_2 + 4a u_1^2 + 4a u_2^2 \big) } },
\]
hence \eqref{eq:S1-1aa_SympInt1} is equal to
\[
    2\pi\sqrt{1 + a} \int_{X_Z}
        \frac{ F(u_1, u_2, u_4,\ldots, u_9) }
        { u_1^{(a - 1)/2} \sqrt{ (u_1 + a u_1^a)\big(a(a + 1)^2 u_1^{2a} + (5a^2 + 2a + 1) u_1^{a + 1} + 4a u_1^a u_2 + 4a u_1^2 + 4a u_2^2 \big) } }
            \, d\vol_{X_Z}.
\]

As in Section~\ref{subsec:S1-1aa_COV2}, it is not feasible to give a complete list of expressions for the Hilbert measure $\Lambda_{X_Z}$, so we choose a pair
of representatives. From the expressions \eqref{eq:S1-1aa_SympRInvarsInW} for the invariants on the fundamental domain, we compute
\begin{align*}
    du_1\wedge du_2 \wedge du_5 \wedge du_6
    &=  - 4 \sqrt{u_1} u_4 u_7
            \, d\mathbf{w}_{123}\wedge d\alpha
    \qquad\text{and}
    \\
    du_1\wedge du_4 \wedge du_8 \wedge du_9
    &=  - \frac{2}{a} \sqrt{u_1} u_5 (u_1 - 2a u_2)
            \, d\mathbf{w}_{123}\wedge d\alpha.
\end{align*}
The expressions for $\Lambda_{X_Z}$ in Equation~\eqref{eq:Main_S1-1aa_Symp} follow,
each defined on a dense subset of the principal stratum, completing the proof of Theorem~\ref{mainthm:S1-1aa_Symp}.


\section{$\OO_m$ acting on $k\R^m$}
\label{sec:Om_kRm}

In this section, we prove the first part of Theorem~\ref{mainthm:Om-SOm}.
Let $\OO_m$ denote the real orthogonal group and let $V_k\simeq\R^{km}$ denote $k$ copies of the defining representation of $\OO_m$.
We denote a point $V_k$ with coordinates $(\mathbf{v}_1,\ldots,\mathbf{v}_k)$ where each $\mathbf{v}_i = (v_{i,1},\ldots,v_{i,m})\in\R^m$ and let $[k] = \{1,2,\ldots,k\}$.
By the First and Second Fundamental Theorems of Invariant Theory for $\OO_m$, see \cite[Sec.~9.2 and 9.4]{PopovVinberg} or \cite[Ch.~11, Sec.~2.1 and 5.2]{ProcesiLieGroups},
a Hilbert basis is given by the collection of $u_{i,j} = \langle\mathbf{v}_i,\mathbf{v}_j\rangle$, $1\leq j\leq i\leq k$, and the relations are generated by
the $(m+1)\times(m+1)$ minors of the Gram matrix $G_k = (u_{i,j})_{i,j\in[k]}$; note that the Hilbert bases working over $\R$ and $\C$ coincide by
\cite[Thrm.~5.8(1)]{GWSliftingHomotopies}. We fix this Hilbert basis as well as the Hilbert embedding
\begin{equation}
\label{eq:Om-HilbertEmbed}
    H\colon V_k\to\R^{{k+1\choose 2}},
    \qquad
    H(\mathbf{v}_1,\ldots,\mathbf{v}_k) = \mathbf{u} = (u_{i,j})_{1\leq j\leq i\leq k},
\end{equation}
where the $u_{i,j}$ are ordered lexicographically. Note that we use $G_k$ to indicate the symmetric matrix $(u_{i,j})_{i,j\in[k]}$ and $\mathbf{u}$ to denote
the element of $\R^{{k+1\choose 2}}$ whose entries are the distinct entries of $G_k$.

By \cite[Thrm~11.18]{GWSlifting}, the representation $V_k$ is \emph{coregular}, meaning that the invariants $u_{i,j}$ are algebraically independent,
if and only if $k\leq m$. In the coregular case, the image $X = H(V_k)$ of the Hilbert embedding is a region in $\R^{{k+1\choose 2}}$ with nonempty interior;
in the non-coregular case, the relations cut down $X$ to be contained in a proper variety of $\R^{{k+1\choose 2}}$. The coregular case
was treated in \cite{HerbigSeatonWhitesell}, so we will henceforth assume that $k > m$ and let $r = k - m \geq 1$.
Section~\ref{subsec:Om_COV1} closely follows that reference.


\subsection{The image of the Hilbert embedding}
\label{subsec:Om_ImageHilb}

The inequalities defining $X$ were determined in \cite[Ex.~0.8]{ProcesiSchwarz}; they are given by declaring that the principal $\ell\times\ell$ minors
of the Gram matrix $G_k$ are nonnegative for $\ell\leq m$. This yields $\sum_{\ell=1}^m {k \choose \ell}$ inequalities, though this number can be reduced as follows.

Let $E_j = E_j(G_k)$ denote the sum of the ${k\choose j}$ principal $j\times j$ minors of $G_k$ and then the relations imply $E_j = 0$ for $j > m$.
By \cite[Eq.~(1.2.13), pp.~53]{HornJohnsonMatrixAnBook}, the characteristic polynomial of $G_k$ is given by
\[
    p_{G_k}(t) = t^k - E_1 t^{k-1} + \cdots + (-1)^m E_m t^r.
\]
Noting that if some $E_j$ vanishes with $j\leq m$, then each subsequent $E_\ell$ must vanish as well by their definition,
we have by \cite[Cor.~7.2.4]{HornJohnsonMatrixAnBook} that $G_k$ is positive definite if and only if each $E_j \geq 0$, $j=1,\ldots,m$.

For subsets $I,J\subseteq[k]$ of the same size $\abs{I} = \abs{J}$,
let $u_{I,J}$ denote the submatrix of $(u_{i,j})$ consisting of rows with index in $I$ and columns with index in $J$.
The image $X$ of the Hilbert embedding \eqref{eq:Om-HilbertEmbed} is given by
\begin{equation}
\label{lem:Om_ImHilb}
    X   =
    \Big\{
        \mathbf{u}\in\R^{{k+1\choose 2}} :
        \det(u_{I,J}) = 0 \;\forall I,J \subseteq[k] \text{ with } \abs{I} = \abs{J} = m + 1,\;
        E_j(G_k) \geq 0 \text{ for } j = 1,2,\ldots, m
    \Big\}.
\end{equation}


\subsection{The first change of variables}
\label{subsec:Om_COV1}

Let $d\mathbf{v} = dv_{1,1}\wedge\cdots\wedge dv_{k,m}$ denote the volume form on $V_k$ where the variables are ordered lexicographically.
We begin with an integral
\begin{equation}
\label{eq:Om_Int0}
    \int_{V_k} f(\mathbf{v}_1, \ldots, \mathbf{v}_k) \, d\mathbf{v}
\end{equation}
where $f$ is a smooth, real-valued, integrable, $\OO_m$-invariant function on $V_k$. We first identify a fundamental domain of a dense subset of $V_k$ and use it
to parametrize an open dense subset of $V_k$. This was carried out in detail in \cite[Sec.~3.1]{HerbigSeatonWhitesell} for the case $k \leq m$, and the reader is referred
to that reference for more details. Here, we will give a summary and explain the modification for the case $k > m$. Note that we assume that the first $m$ vectors
$\mathbf{v}_1, \ldots, \mathbf{v}_m$ are linearly independent, excluding a set of measure zero.

The fundamental domain is given by
\begin{equation}
\label{eq:Om_FundDomain}
    W   =   \big\{(w_{i,j}) = (\mathbf{w}_1,\ldots,\mathbf{w}_k)\in\R^{km} :
                w_{i,j} = 0 \text{ if } j > i, \;
                w_{i,i} > 0 \text{ for } 1\leq i\leq m \big\}.
\end{equation}
Each $(\mathbf{v}_1,\ldots,\mathbf{v}_k)\in V_k$ can be transported to an element of $W$ through the following iterative process.
Let $R_\theta^{j}\in\SO_m$ be the rotation through $\theta$ acting on coordinates $j$ and $j+1$ of $\R^m$, leaving all other coordinates fixed.
We first apply a rotation $R_{-\theta^{(1,1)}}^{(m-1)}$, where the angle $\theta^{(1,1)}\in(-\pi,\pi]$ is chosen such that
$R_{-\theta^{(1,1)}}^{(m-1)}\mathbf{v}_1$ has vanishing $m$\textsuperscript{th} coordinate and nonnegative $(m-1)$\textsuperscript{st} coordinate.
We then apply rotation $R_{-\theta^{(1,2)}}^{(m-2)}$ with $\theta^{(1,2)}\in[0,\pi)$ so that
$R_{-\theta^{(1,2)}}^{(m-2)} R_{-\theta^{(1,1)}}^{(m-1)}\mathbf{v}_1$ has vanishing $m$\textsuperscript{th} and $(m-1)$\textsuperscript{st}
coordinates and nonnegative $(m-2)$\textsuperscript{nd} coordinate. Continuing, we apply a total $m-1$ rotations to obtain rotated vectors
$( \mathbf{v}_1^{(1)},\ldots,\mathbf{v}_k^{(1)} )$ such
that,\footnote{After the $j$\textsuperscript{th} rotation, it need not be the case that the
$(m-j-1)$\textsuperscript{st} coordinate of the transformed vector is nonzero as was incorrectly stated in \cite{HerbigSeatonWhitesell}.
However, the linear independence assumption ensures that $v_{1,1}^{(1)}$ will be nonzero after the first sequence of $m-1$ rotations.}
letting $v_{1,i}^{(1)}$ denote the coordinates of $\mathbf{v}_1^{(1)}$ as above, $v_{1,1}^{(1)} > 0$ and $v_{1,i}^{(1)} = 0$ for $i > 1$.
The corresponding transformation of the volume form \cite[pp.~7]{HerbigSeatonWhitesell} is
\[
    d\mathbf{v}
    =
    (v_{1,1}^{(1)})^{m-1} \prod\limits_{j=2}^{m-1} \sin^{j-1}\theta^{(1,j)}
        \, dv_{1,1}^{(1)} \wedge \bigwedge\limits_{j=1}^{m-1} d\theta^{(1,m-j)} \wedge \bigwedge\limits_{\substack{2\leq i\leq k \\ 1\leq j\leq m}} dv_{i,j}^{(1)},
\]

Next, we apply a similar sequence of $m-2$ rotations through angles $\theta^{(2,1)},\ldots,\theta^{(2,m-2)}$,
each chosen so that a coordinate of the second vector vanishes beginning with the last,
and each fixing the first coordinate so that $\mathbf{v}_1^{(1)}$ is fixed. The resulting element $( \mathbf{v}_1^{(2)},\ldots,\mathbf{v}_k^{(2)} )$ of $V_k$ satisfies
$v_{2,2}^{(2)} > 0$ and $v_{(2,i)}^{(2)} = 0$ for $i > 2$. We apply a similar process targeting $\mathbf{v}_3$, etc., resulting in 
$( \mathbf{v}_1^{(m-1)},\ldots,\mathbf{v}_k^{(m-1)} ) \in V_k$. Finally, we apply a reflection if necessary to ensure that the last coordinate
of the $m$\textsuperscript{th} vector is positive and denote the resulting element $(\mathbf{w}_1,\ldots,\mathbf{w}_k)\in V_k$.
Note that the process here coincides with that described in \cite{HerbigSeatonWhitesell} for $k = m$ and then terminates.
This process yields an element $R\in\OO_m$ such that multiplication by $R$ yields a map
\[
    \begin{pmatrix}
        v_{1,1} &   v_{2,1} &   \cdots  &   v_{m,1} &   \cdots  &   v_{m+r,1}   \\
        v_{1,2} &   v_{2,2} &   \cdots  &   v_{m,2} &   \cdots  &   v_{m+r,2}   \\
        \vdots  &   \vdots  &   \ddots  &   \vdots  &   \cdots  &   \ddots      \\
        v_{1,m} &   v_{2,m} &   \cdots  &   v_{m,m} &   \cdots  &   v_{m+r,m}
    \end{pmatrix}
    \overset{R}{\mapsto}
    \begin{pmatrix}
        w_{1,1} &   w_{2,1} &   \cdots  &   w_{m,1} &   \cdots  &   w_{m+r,1}   \\
        0       &   w_{2,2} &   \cdots  &   w_{m,2} &   \cdots  &   w_{m+r,2}   \\
        \vdots  &   \vdots  &   \ddots  &   \vdots  &   \cdots  &   \ddots      \\
        0       &   0       &   \cdots  &   w_{m,m} &   \cdots  &   w_{m+r,m}
    \end{pmatrix}
    \in W,
\]
and the transformation of the volume form is computed in \cite[pp.~7]{HerbigSeatonWhitesell} to be
\[
    d\mathbf{v}
    =
    2 \prod\limits_{i=1}^{m-1} w_{i,i}^{m-i}
    \prod\limits_{i=1}^{m-1} \prod\limits_{j=1}^{m-i} \sin^{j-1}\theta^{(i,j)}
        \bigwedge\limits_{i=1}^k \left( \bigwedge\limits_{j=1}^{\min\{i,m\}} dw_{i,j}
            \wedge \bigwedge\limits_{\ell=1}^{m-1} d\theta^{(i,m-\ell)}
        \right).
\]
Setting
\[
    \Theta  =
    \Big\{ \big(\theta^{(i,j)}\big)_{\begin{subarray}{l}1\leq i\leq m-1 \\ 1\leq j\leq m-i\end{subarray}} \in (-\pi,\pi)^{(m^2-m)/2}: 
        (\theta^{(i,j)}) \in [0,\pi) \text{ for } j > 1\Big\}
\]
and
\[
    \mathbf{w} = (w_{ij})_{\begin{subarray}{l}i\in[m] \\ 1\leq j \leq \min\{i,m\}\end{subarray}}\in\R^{{m+1\choose 2}+rm}
\]
ordered lexicographically, we may rewrite the integral \eqref{eq:Om_Int0} as
\begin{align*}
    &2\int_{W\times\Theta}  f(\mathbf{w})
        \prod\limits_{i=1}^{m-1} w_{i,i}^{m-i}
        \prod\limits_{i=1}^{m-1} \prod\limits_{j=1}^{m-i} \sin^{j-1}\theta^{(i,j)}
        \bigwedge\limits_{i=1}^k \left( \bigwedge\limits_{j=1}^{\min\{i,m\}} dw_{i,j}
            \wedge \bigwedge\limits_{\ell=1}^{m-1} d\theta^{(i,m-\ell)}
        \right)
    \\&\qquad =
    2\int_{W}  f(\mathbf{w})
        \prod\limits_{i=1}^{m-1} w_{i,i}^{m-i} \,
        \bigwedge\limits_{\substack{1\leq i\leq k \\ 1\leq j \leq {\min\{i,m\}} }}  dw_{i,j}
    \int_\Theta \prod\limits_{i=1}^{m-1} \prod\limits_{j=1}^{m-i} \sin^{j-1}\theta^{(i,j)}
        \, \bigwedge\limits_{\ell=1}^{m-1} d\theta^{(i,m-\ell)}.
\end{align*}
To compute the integral over $\Theta$, we have by \cite[6.2.1 and 6.2.2]{AbramowitzStegun} that
\[
    \int\limits_0^\pi \sin^\ell \theta\;d\theta
        =   \frac{\sqrt{\pi}\;\Gamma\left(\frac{\ell+1}{2}\right)}{\Gamma\left(\frac{\ell+2}{2}\right)},
\]
and then
\begin{align*}
    2\int_\Theta \prod\limits_{i=1}^{m-1} \prod\limits_{j=1}^{m-i} \sin^{j-1}\theta^{(i,j)}
        \, \bigwedge\limits_{\ell=1}^{m-1} d\theta^{(i,m-\ell)}
    &=
    2^m\pi^{(m-1)(m+2)/4}
        \prod\limits_{i=1}^{m-2} \left( \frac{ \Gamma\left(\frac{i+1}{2}\right) }{ \Gamma\left(\frac{i+2}{2}\right)} \right)^{m-i-1}
    \\&=    \prod\limits_{i=1}^m \Vol(S^{i-1})
    \qquad= 2^m \gamma
\end{align*}
where $\Vol(S^{\ell-1}) = \frac{2\pi^{\ell/2}}{\Gamma(\ell/2)}$ is the volume of the $(\ell-1)$-dimensional unit sphere $S^{\ell-1}$ and
\[
    \gamma = \pi^{(m-1)(m+2)/4}
        \prod\limits_{i=1}^{m-2} \left( \frac{ \Gamma\left(\frac{i+1}{2}\right) }{ \Gamma\left(\frac{i+2}{2}\right)} \right)^{m-i-1}.
\]
Then \eqref{eq:Om_Int0} is equal to
\begin{equation}
\label{eq:Om_Int1}
    \prod\limits_{i=1}^m \Vol(S^{i-1}) \int_{W}  g(\mathbf{w})
        \prod\limits_{i=1}^{m-1} w_{i,i}^{m-i} \,
        \bigwedge\limits_{\substack{1\leq i\leq k \\ 1\leq j \leq {\min\{i,m\}} }}  dw_{i,j},
\end{equation}
where $g(\mathbf{w}) = f(\mathbf{v})$ and the $dw_{i,j}$ are ordered lexicographically.


\subsection{The second change of variables}
\label{subsec:Om_COV2}

There does not appear to be a simple closed form expression for $d\vol_X$, so we turn to the computation of the Hilbert measure.
Specifically, we express $\bigwedge_{\begin{subarray}{l}1\leq i\leq k \\ 1\leq j \leq {\min\{i,m\}} \end{subarray}}  dw_{i,j}$ in terms of the
$u_{i,j}$ coordinates. Throughout this section, we assume the $(\mathbf{w}_1,\ldots,\mathbf{w}_k)\in V_k$ is contained in the fundamental
domain $W$ described in \eqref{eq:Om_FundDomain}.

Let $G_\ell = (\langle \mathbf{w}_i, \mathbf{w}_j\rangle)_{i,j\in[\ell]}$ denote the Gram matrix of the first $\ell$ vectors
$\mathbf{w}_1,\ldots,\mathbf{w}_\ell$. Note that this is consistent with the notation $G_k$ and $G_m$ above, as $G_\ell$
is as well equal to $(\langle \mathbf{v}_i, \mathbf{v}_j\rangle)_{i,j\in[\ell]}$ as the scalar products are $\OO_m$-invariants.
As the span of the $\mathbf{w}_i$ has dimension $m$, $\det G_\ell = 0$ for $\ell > m$. The proof of
\cite[Lem.~3.1]{HerbigSeatonWhitesell} adapts readily to the present case, which we include for the benefit of the reader.

\begin{lemma}
\label{lem:Om_COV2_wii}
With $\mathbf{w}_i$ and $G_\ell$ as above, we have
\begin{equation}
\label{eq:Om_COV2_wii}
    w_{i,i} =   \sqrt{ \frac{\det G_i}{\det G_{i-1}} }
\end{equation}
for $i = 1, 2, \ldots, m$, where we set $\det G_0 = 1$.
\end{lemma}
\begin{proof}
The proof is by induction on $i$. If $i = 1$, then $w_{1,1}^2 = \langle\mathbf{w}_1,\mathbf{w}_1\rangle = \det G_1$ so that
Equation~\eqref{eq:Om_COV2_wii} holds. So assume that for some $i \leq m$, Equation~\eqref{eq:Om_COV2_wii} holds for each
$w_{j,j}$ with $1\leq j < i$. Let $W_i$ be the $i\times i$ matrix
\begin{equation}
\label{eq:Om_COV2_wiiMatrix}
    W_i =
    \begin{pmatrix}
        w_{1,1} &   w_{2,1} &   \cdots  &   w_{i,1}   \\
        0       &   w_{2,2} &   \cdots  &   w_{i,2}   \\
        \vdots  &   \vdots  &   \ddots  &   \vdots    \\
        0       &   0       &   \cdots  &   w_{i,i}
    \end{pmatrix}
\end{equation}
and then $G_i = W_i^T W_i$ so that by the inductive hypothesis,
\begin{align*}
    \det G_i
    &=      (\det W_i)^2
    \\&=    \prod\limits_{j=1}^i w_{j,j}^2
    \\&=    w_{i,i}^2 \prod\limits_{j=1}^{i-1} \frac{\det G_j}{\det G_{j-1}}
    \\&=    w_{i,i}^2 \det G_{i-1}.
    \qedhere
\end{align*}
\end{proof}

We now have the following.

\begin{lemma}
\label{lem:Om_COV2_SolveUforW}
An element $(\mathbf{w}_1,\ldots,\mathbf{w}_k)\in W$ is determined by the
$u_{i,j} = \langle\mathbf{w}_i,\mathbf{w}_j\rangle$ with $i\in[k]$ and $1\leq j \leq \min\{i,m\}$.
\end{lemma}
\begin{proof}
The $u_{i,j}$ with $1\leq j \leq i\leq m$ correspond to the coordinates of the Hilbert embedding \eqref{eq:Om-HilbertEmbed} in the case $k = m$,
i.e., ignoring the vectors $\mathbf{w}_{m+1}, \ldots, \mathbf{w}_k$. It follows that these coordinates determine the $\OO_m$-orbit of
$(\mathbf{w}_1,\ldots,\mathbf{w}_m)$. Then as the $\OO_m$-orbit contains only one element in $W$, the vectors $\mathbf{w}_1,\ldots,\mathbf{w}_m$
are determined. 

Now, given $\mathbf{w}_1,\ldots,\mathbf{w}_m$, the definitions of the $u_{i,j}$ with $j\in[m]$ and $m + 1 \leq i \leq k$, treated as a linear system
whose indeterminates are the coordinates of $\mathbf{w}_{m+1},\ldots,\mathbf{w}_k$, is triangular with each diagonal entry $w_{i,i} > 0$ and hence
can be solved by back-substitution.
\end{proof}

As a consistency check, note that by the Corollary of \cite[Thrm.~1]{Kutz}, the ideal in the (regular, hence Cohen-Macaulay) ring $S$ of symmetric
$m\times m$ matrices generated by the $(m+1)\times(m+1)$ minors has height ${k - m + 1 \choose 2}$, so that the corresponding variety has dimension
\[
    {k + 1 \choose 2} - {k - m + 1 \choose 2}
    =   \frac{m(2k - m + 1)}{2}.
\]
This coincides with the number of $u_{ij}$ identified in Lemma~\ref{lem:Om_COV2_SolveUforW}.

It follows that, letting
\begin{equation}
\label{eq:Om_COV2_DefUHat}
    \hat{\mathbf{u}} = (u_{i,j})_{\begin{subarray}{l} i\in[k] \\ 1\leq j \leq \min\{i,m\}\end{subarray}}
\end{equation}
ordered lexicographically, we may express
\begin{equation}
\label{eq:Om_COV2_dwTodu}
    \det\Big(\frac{\partial \hat{\mathbf{u}} }{\partial \mathbf{w} }\Big)
        \bigwedge\limits_{\substack{1\leq i\leq k \\ 1\leq j \leq {\min\{i,m\}} }}  dw_{i,j}
    =
    \bigwedge\limits_{\substack{1\leq i\leq k \\ 1\leq j \leq {\min\{i,m\}} }}  du_{i,j}.
\end{equation}

It remains only to compute $\det\big(\frac{\partial \hat{\mathbf{u}} }{\partial \mathbf{w} }\big)$.
Once again, this computation is a minor modification of \cite[Lem.~3.2]{HerbigSeatonWhitesell}.
Recall that $r = k - m$.

\begin{lemma}
\label{lem:Om_COV2_JacobianDet}
The Jacobian $\frac{\partial \hat{\mathbf{u}} }{\partial \mathbf{w} }$ is lower-triangular, and
\[
    \det\Big(\frac{\partial \hat{\mathbf{u}} }{\partial \mathbf{w} }\Big)
    =
    2^m \sqrt{\det G_1\det G_2\cdots\det G_{m-1}(\det G_m)^{r+1} }.
\]
\end{lemma}
\begin{proof}
Note that the indices of the elements of $\mathbf{w}$ and $\hat{\mathbf{u}}$ coincide, hence, the diagonal entries of
$\frac{\partial \hat{\mathbf{u}} }{\partial \mathbf{w} }$ are given by $\frac{\partial u_{i,j} }{\partial w_{i,j} }$.
When $j\leq \min\{i,m\}$, $u_{i,j} = w_{i,1} w_{j,1} + \cdots + w_{i,j} w_{j,j}$ depends only on $w_{i,s}$ and $w_{j,s}$ with $s\leq j$.
It follows that $\frac{\partial \hat{\mathbf{u}} }{\partial \mathbf{w} }$ is
lower-triangular. By Lemma~\ref{lem:Om_COV2_wii}, the diagonal entries are
\begin{equation}
\label{eq:Om_COV2_JacobianDetPartialDeriv}
    \frac{\partial u_{i,j}}{\partial w_{i,j}}
    =
    (1 + \delta_{i,j}) w_{j,j}
    =
    (1 + \delta_{i,j}) \sqrt{ \frac{\det G_j}{\det G_{j-1}} }.
\end{equation}
where $\delta_{i,j}$ is the Kronecker delta. Hence
\begin{align*}
    \det\Big(\frac{\partial \hat{\mathbf{u}} }{\partial \mathbf{w} }\Big)
    &=
    \prod\limits_{i=1}^k \prod\limits_{j=1}^{\min\{i,m\}} (1 + \delta_{i,j}) \sqrt{ \frac{\det G_j}{\det G_{j-1}} }
    \\&=
    2^m \prod\limits_{i=1}^k \prod\limits_{j=1}^{\min\{i,m\}} \sqrt{ \frac{\det G_j}{\det G_{j-1}} }
    \\&=
    2^m \prod\limits_{i=1}^k  \sqrt{ \det G_{\min\{i,m\}} }
    \\&=
    2^m \sqrt{\det G_1\det G_2\cdots\det G_{m-1}(\det G_m)^{r+1} }.
    \qedhere
\end{align*}
\end{proof}

Returning to the integral in \eqref{eq:Om_Int1}, we use Equation~\eqref{eq:Om_COV2_dwTodu} and Lemmas~\ref{lem:Om_COV2_wii} and \ref{lem:Om_COV2_JacobianDet}
to compute
\begin{align*}
    \prod\limits_{i=1}^{m-1} w_{i,i}^{m-i} \,
        \bigwedge\limits_{\substack{1\leq i\leq k \\ 1\leq j \leq {\min\{i,m\}} }}  dw_{i,j}
    &=
    \frac{\prod\limits_{i=1}^{m-1} \frac{(\det G_i)^{(m - i)/2}}{(\det G_{i-1})^{(m - i)/2} } }
        { 2^m \sqrt{\det G_1\det G_2\cdots\det G_{m-1}(\det G_m)^{r+1} } }
        \bigwedge\limits_{\substack{1\leq i\leq k \\ 1\leq j \leq {\min\{i,m\}} }}  du_{i,j}
    \\&=
    \frac{\sqrt{\det G_1\det G_2\cdots\det G_{m-1} } }
        { 2^m \sqrt{\det G_1\det G_2\cdots\det G_{m-1}(\det G_m)^{r+1} } }
        \bigwedge\limits_{\substack{1\leq i\leq k \\ 1\leq j \leq {\min\{i,m\}} }}  du_{i,j}
    \\&=
    \frac{1}
        { 2^m (\det G_m)^{(r+1)/2} }
        \bigwedge\limits_{\substack{1\leq i\leq k \\ 1\leq j \leq {\min\{i,m\}} }}  du_{i,j}.
\end{align*}
Applying \cite[Thrm.~1]{GWSdiff} to obtain a smooth function $F$ on $\R^{{k+1\choose 2}}$ such that $F(\mathbf{u}) = g(\mathbf{w}) = f(\mathbf{v})$,
\eqref{eq:Om_Int1} can be written in the form
\begin{equation}
\label{eq:Om_Int2}
    \frac{1}{2^m}\prod\limits_{i=1}^m \Vol(S^{i-1}) \int_{X}  \frac{ F(\mathbf{u}) }{ (\det G_m)^{(r+1)/2}}
         \, \bigwedge\limits_{\substack{1\leq i\leq k \\ 1\leq j \leq {\min\{i,m\}} }}  du_{i,j}.
\end{equation}
It follows that the Hilbert measure is that given in Equation~\eqref{eq:Main_Om}, completing the proof of the first part of Theorem~\ref{mainthm:Om-SOm}.


\section{$\SO_m$ acting on $k\R^m$}
\label{sec:SOm_kRm}

We now turn to the second part of Theorem~\ref{mainthm:Om-SOm}.
Let $\SO_m$ denote the real special orthogonal group and $V_{\text{SO},k}\simeq\R^{km}$ denote $k$ copies of the defining representation of $\SO_m$, i.e.,
$V_{\text{SO},k}$ is the restriction of the representation $V_k$ in Section~\ref{sec:Om_kRm} to the subgroup $\SO_m \leq \OO_m$. We use the same coordinates
and notation as in Section~\ref{sec:Om_kRm}. The computation of the Hilbert measure in this case is very similar to that of the $\OO_m$ representation $V_k$,
so we focus on the differences.

By the First and Second Fundamental Theorems of Invariant Theory for $\SO_m$, see \cite[Sec.~9.2 and 9.4]{PopovVinberg} or \cite[Ch.~11, Sec.~2.1 and 5.2]{ProcesiLieGroups},
a Hilbert basis is given by the $\OO_m$-invariants $u_{i,j} = \langle\mathbf{v}_i,\mathbf{v}_j\rangle$, $1\leq j\leq i\leq k$ as well as the determinants
$\det\nolimits_{i_1,\ldots,i_m} = \det(\mathbf{v}_{i_1},\ldots,\mathbf{v}_{i_m})$, $1\leq i_1 < \cdots < i_m \leq k$, and the relations are generated by
the $(m+1)\times(m+1)$ minors of the Gram matrix $G_k = (u_{i,j})_{i,j\in[k]}$ as well as relations of the following form. If $1\leq i_1 < \cdots < i_m \leq k$
and $1\leq j_1 < \cdots < j_m \leq k$, then
\[
    \det( u_{i_a, j_b} )_{a,b \in [k]} = \det\nolimits_{i_1,\ldots,i_m} \det\nolimits_{j_1,\ldots,j_m}.
\]
We fix this Hilbert basis as well as the Hilbert embedding
\begin{equation}
\begin{split}
\label{eq:SOm-HilbertEmbed}
    &H \colon V_k\to\R^{{k+1\choose 2}+{k\choose m}},
    \\
    &H(\mathbf{v}_1,\ldots,\mathbf{v}_k) = (\mathbf{u},\mathbf{\det}) = \big((u_{i,j})_{1\leq j\leq i\leq k},
        (\det(\mathbf{v}_{i_1},\ldots,\mathbf{v}_{i_m}))_{1\leq i_1 < \cdots < i_m \leq k}\big).
\end{split}
\end{equation}
Again by \cite[Thrm~11.18]{GWSlifting}, $V_k$ is coregular if and only if $k\leq m - 1$, and the coregular case was treated in \cite{HerbigSeatonWhitesell}.
We therefore assume that $k > m - 1$. We therefore always have at least one invariant of the form $\det\nolimits_{i_1,\ldots,i_m}$. Note that the invariant
$\det\nolimits_{[m]}$ will play an important role in these computations, as will the relation
\begin{equation}
\label{eq:SOm-MainRelation}
    \det G_m = \det\nolimits_{[m]}^2.
\end{equation}

The inequalities defining the image $X$ of the Hilbert embedding \eqref{eq:SOm-HilbertEmbed} include the inequalities from Section~\ref{subsec:Om_ImageHilb}
as well as new inequalities involving the scalar products of gradients of the $\det\nolimits_{i_1,\ldots,i_m}$ invariants with one another as well as with the $u_{i,j}$.
They can be computed in examples using \cite[Thrm~0.10 and Prop.~0.6]{ProcesiSchwarz}, and the list of inequalities can be shortened using
\cite[Cor.~7.2.4]{HornJohnsonMatrixAnBook} as in Section~\ref{subsec:Om_ImageHilb}; they do not seem conducive to being written down explicitly in general.

Beginning with an integral of the form in \eqref{eq:Om_Int0}, where $f$ is now a smooth, real-valued, integrable, $\SO_m$-invariant function on $V_k$,
the first change of variables is almost identical to that in Section~\ref{subsec:Om_COV1}. The only change occurs at the last step, where a reflection
was used to ensure that the $m$\textsuperscript{th} coordinate of the $m$\textsuperscript{th} vector is positive. This is no longer possible, so the fundamental
domain becomes
\begin{equation}
\label{eq:SOm_FundDomain}
    W   =   \big\{(w_{i,j}) = (\mathbf{w}_1,\ldots,\mathbf{w}_k)\in\R^{km} :
                w_{i,j} = 0 \text{ if } j > i, \;
                w_{i,i} > 0 \text{ for } 1\leq i\leq m - 1, w_{m,m}\neq 0 \big\}.
\end{equation}
Transforming the integral to one over the fundamental domain, this changes only the final factor of $2 = \Vol(S^0)$, so the resulting integral is in this case
\begin{equation}
\label{eq:SOm_Int1}
    \prod\limits_{i=2}^m \Vol(S^{i-1}) \int_{W}  f(\mathbf{w})
        \prod\limits_{i=1}^{m-1} w_{i,i}^{m-i} \,
        \bigwedge\limits_{\substack{1\leq i\leq k \\ 1\leq j \leq {\min\{i,m\}} }}  dw_{i,j}.
\end{equation}

As in Section~\ref{subsec:Om_COV2}, let $(\mathbf{w}_1,\ldots,\mathbf{w}_k)\in W\subset V_k$ and let
$G_\ell = (\langle \mathbf{w}_i, \mathbf{w}_j\rangle)_{i,j\in[\ell]}$ denote the Gram matrix of the first $\ell$ vectors.
In this case, the proof of Lemma~\ref{lem:Om_COV2_wii} carries through identically except for the $m$\textsuperscript{th} step,
where the fact that $w_{m,m}$ may be negative introduces a sign. Using the definition of $W_m$ in Equation~\eqref{eq:Om_COV2_wiiMatrix},
this sign is seen to be the sign of the invariant $\det\nolimits_{[m]}$, so keeping track of this sign in the proof of
Lemma~\ref{lem:Om_COV2_wii} yields the following.

\begin{lemma}
\label{lem:SOm_COV2_wii}
With $\mathbf{w}_i$ and $G_\ell$ as above, we have
\begin{equation}
\label{eq:SOm_COV2_wii}
    w_{i,i} =   \sqrt{ \frac{\det G_i}{\det G_{i-1}} }
    \text{ for } 1\leq i\leq m-1,
    \qquad\text{ and }\qquad
    w_{m,m} =   \sign(\det\nolimits_{[m]})\sqrt{ \frac{\det G_m}{\det G_{m-1}} },
\end{equation}
where we set $\det G_0 = 1$.
\end{lemma}

We have the following minor variation of Lemma~\ref{lem:Om_COV2_SolveUforW}.

\begin{lemma}
\label{lem:SOm_COV2_SolveUforW}
An element $(\mathbf{w}_1,\ldots,\mathbf{w}_k)\in W$ is determined by the sign of $\det\nolimits_{[m]}$ and the
$u_{i,j} = \langle\mathbf{w}_i,\mathbf{w}_j\rangle$ with $i\in[k]$ and $1\leq j \leq \min\{i,m\}$.
\end{lemma}
\begin{proof}
In the case $k = m$, the Hilbert embedding \eqref{eq:SOm-HilbertEmbed} is now given by
the $u_{i,j}$ with $1\leq j \leq i\leq m$ as well as $\det\nolimits_{[m]}$, so as in the proof of Lemma~\ref{lem:Om_COV2_SolveUforW},
these values determine the $\SO_m$-orbits of the first $m$ vectors. Using the relation in Equation~\eqref{eq:SOm-MainRelation},
$\abs{\det\nolimits_{[m]}}$ is determined by the $u_{i,j}$ with $1\leq j \leq i\leq m$, so only the sign of $\det\nolimits_{[m]}$
is required. The rest of the proof proceeds as in the proof of Lemma~\ref{lem:Om_COV2_SolveUforW}.
\end{proof}

We therefore continue to let $\hat{\mathbf{u}}$ be as defined in Equation~\eqref{eq:Om_COV2_DefUHat}. Then
the only change to Lemma~\ref{lem:Om_COV2_JacobianDet} is consideration of the sign of $\det\nolimits_{[m]}$ when applying
Lemma~\ref{lem:SOm_COV2_wii}. Recall that $r = k - m$.

\begin{lemma}
\label{lem:SOm_COV2_JacobianDet}
The Jacobian $\frac{\partial \hat{\mathbf{u}} }{\partial \mathbf{w} }$ is lower-triangular, and
\[
    \det\Big(\frac{\partial \hat{\mathbf{u}} }{\partial \mathbf{w} }\Big)
    =
    \sign(\det\nolimits_{[m]})^{r+1} 2^m \sqrt{\det G_1\det G_2\cdots\det G_{m-1}(\det G_m)^{r+1} }.
\]
\end{lemma}
\begin{proof}
Following the proof of Lemma~\ref{lem:Om_COV2_JacobianDet}, Equation~\eqref{eq:Om_COV2_JacobianDetPartialDeriv} is unchanged when $(i, j)\neq (m, m)$
and becomes
\[
    \frac{\partial u_{m,m}}{\partial w_{m,m}}
    =
    2 w_{m,m}
    =
    2 \sign(\det\nolimits_{[m]})\sqrt{ \frac{\det G_m}{\det G_{m-1}} }
\]
in this case. In the computation of $\det( \partial \hat{\mathbf{u}} / \partial \mathbf{w})$, this factor occurs $r + 1$ times.
\end{proof}

With this, a simple inclusion of $\sign(\det\nolimits_{[m]})^{r+1}$ throughout the denominators in the computation of $\Lambda$ at the end of
Section~\ref{subsec:Om_COV2} yields the expression for $\Lambda_X$ in Equation~\eqref{eq:Main_SOm} and completes the proof of Theorem~\ref{mainthm:Om-SOm}.


\section*{AI Disclosure}

The contents of this paper were human-generated. The third author used ChatGPT (OpenAI) and Claude (Anthropic) in the late stages of manuscript preparation
to identify references, fact-check and organize background information, and check completed drafts for correctness and typographical errors.
In particular, these tools identified minor errors and suggested the simplified proof of Lemma~\ref{lem:S1-1aa_ImHilb} as well as the
technique used to reduce the number of inequalities in Section~\ref{subsec:Om_ImageHilb}. After using these tools, the authors
addressed corrections manually and take full responsibility for the contents.


\bibliographystyle{amsplain}
\bibliography{LRS-PS}

\end{document}